\documentclass[11pt,a4paper]{article}

\usepackage{amssymb,amsthm,amsmath}
\usepackage[pdftex]{graphicx}
\usepackage{graphicx}
\usepackage[margin=1in]{geometry}
\usepackage{color}
\usepackage{hyperref}
\usepackage{comment}
\usepackage{bm}
\usepackage{cleveref}

\newtheorem{theorem}{Theorem}[section]
\newtheorem{lemma}[theorem]{Lemma}
\newtheorem{proposition}[theorem]{Proposition}

\theoremstyle{definition}
\newtheorem*{definition}{Definition}

\theoremstyle{remark}

\newcommand{\cB}{\mathcal{B}}

\newcommand{\cK}{\mathcal{K}}

\newcommand{\cM}{\mathcal{M}}

\newcommand{\cP}{\mathcal{P}}
\newcommand{\cQ}{\mathcal{Q}}

\newcommand{\cS}{\mathcal{S}}
\newcommand{\cT}{\mathcal{T}}

\newcommand{\cX}{\mathcal{X}}
\newcommand{\cY}{\mathcal{Y}}
\newcommand{\cZ}{\mathcal{Z}}
\newcommand{\CX}{\mathcal{X}}
\newcommand{\CY}{\mathcal{Y}}
\newcommand{\CZ}{\mathcal{Z}}

\newcommand{\bM}{\mathbf{M}}

\newcommand{\bP}{\mathbf{P}}
\newcommand{\bQ}{\mathbf{Q}}
\newcommand{\bR}{\mathbf{R}}
\newcommand{\bS}{\mathbf{S}}

\newcommand{\balpha}{\bm{\alpha}}

\newcommand{\bu}{\mathbf{u}}
\newcommand{\bv}{\mathbf{v}}

\newcommand{\bx}{\mathbf{x}}
\newcommand{\by}{\mathbf{y}}
\newcommand{\bz}{\mathbf{z}}

\newcommand{\bbE}{\mathbb{E}}

\newcommand{\bbX}{\mathbb{X}}
\newcommand{\bbY}{\mathbb{Y}}
\newcommand{\bbZ}{\mathbb{Z}}

\newcommand{\fmin}{f_{\min}}
\newcommand{\RR}{\mathbb{R}}
\newcommand{\NN}{\mathbb{N}}
\newcommand{\iy}{\mathit{y}}
\newcommand{\biy}{\bm{\mathit{\iy}}}
\newcommand{\bbeta}{\bm{\beta}}
\newcommand{\bgamma}{\bm{\gamma}}
\newcommand{\mlb}{\mathrm{mlb}}

\newcommand{\btau}{\bm{\tau}}
\newcommand{\Py}{P}

\newcommand{\FX}{\mathfrak{X}}

\newcommand{\gw}{\mathrm{GW}}
\newcommand{\gwp}{\mathrm{gw}}

\title{Moment Sum-of-Squares Hierarchy for Gromov Wasserstein: Continuous Extensions and Sample Complexity}

\author{Hoang Anh Tran, Binh Tuan Nguyen, Yong Sheng Soh}

\begin{document}

\maketitle

\begin{abstract}
The Gromov-Wasserstein (GW) problem is an extension of the classical optimal transport problem to settings where the source and target distributions reside in incomparable spaces, and for which a cost function that attributes the price of moving resources is not available.  The sum-of-squares (SOS) hierarchy is a principled method for deriving tractable semidefinite relaxations to generic polynomial optimization problems.  In this work, we apply ideas from the moment-SOS hierarchy to solve the GW problem.

More precisely, we identify extensions of the moment-SOS hierarchy, previously introduced for the discretized GW problem in \cite{tran:25}, such that they remain valid for general probability distributions.  This process requires a suitable generalization of positive semidefiniteness over finite-dimensional vector spaces to the space of probability distributions.  We prove the following properties concerning these continuous extensions:  First, these relaxations form a genuine hierarchy in that the optimal value converges to the GW distance.  Second, each of these relaxations induces a pseudo-metric over the collection of metric measure spaces.  Crucially, unlike the GW problem, these induced instances are tractable to compute -- the discrete analogs are expressible as semidefinite programs and hence are tractable to solve.  Separately from these properties, we also establish a statistical consistency result arising from sampling the source and target distributions.  Our work suggests fascinating applications of the SOS hierarchy to optimization problems over probability distributions in settings where the objective and constraint depend on these distributions in a polynomial way.
\end{abstract}

\section{Introduction}

The optimal transport problem (OT) 
is a computational task that concerns moving resources from one location to another in a fashion that minimizes the overall transportation cost.  It finds applications in a wide range of scientific and engineering disciplines.  For instance, the OT problem has direct applications to operations research and economics in which the resource represents a commodity that needs to be transported \cite{galichon2016optimal}.  More recently, the OT problem has found applications in statistics and machine learning in which the resources may be probability distributions (of datasets) or physical objects such as natural images \cite{cuturi2013sinkhorn,kusner2015word,solomon2016entropic,peyre2019computational}.

To set up the optimal transport problem in the classical setting, one typically needs a suitable notion of cost that quantifies the price of moving a unit amount of resource from one location to another.  Subsequently, the key decision variable of interest is a {\em transportation plan} that provides instructions about how the resource should be transported -- for instance, for every unit amount of resource at any location, the plan specifies the proportion as well as the destination to which this resource is to be transported to. 

The Gromov-Wasserstein (GW) problem is an extension of the classical optimal transport problem that applies to settings in which a standard cost function that describes the price of moving resources is unavailable  \cite{Mem:07,Mem:11}.  This can happen in practical situations in which the source and the destination resides in two different spaces that are not naturally compatible with each other.  To circumvent such difficulties, the Gromov-Wasserstein problem provides an alternative formulation that yields reasonably-looking transportation plans.  It first supposes that we have access to the {\em metric} evaluations of the spaces where the source and the target distributions (resources) reside, and second formulates an objective that penalizes for discrepancy between these metric spaces. The Gromov-Wasserstein distance, together with its entropy-regularized version \cite{solomon2016entropic,peyre2016gromov,RGK:24} as well as other related variants \cite{SPC:22,TTZCFDB:22,PSW:25}, has been widely used in numerous computational applications, such as graph learning \cite{xu2019gromov,vayer2020fused,vincent-cuaz_online_2021}, shape matching \cite{Mem:11}, single-cell alignment \cite{demetci2020gromov}, alignment of language embedding vectors \cite{alvarez2018gromov,grave2019unsupervised}, and generative modeling \cite{bunne2019learning}.  The work in \cite{Sturm:06,Sturm:23} studies geometric aspects of the GW distance, while a broader treatment of the topic can be found in \cite{vayer2020contribution,tran2024optimal}.

More concretely, a {\em metric measure space} is a triplet $(\cX,d_{\cX},\mu)$ whereby $\cX$ specifies a set (a {\em space}, so to speak), $d_{\cX}$ specifies a metric over $\cX$, while $\mu$ specifies a probability distribution over $\cX$.  Given a pair of metric measure spaces $\mathbb{X} := (\cX, d_{\cX}, \mu)$ and $\mathbb{Y} := (\cY, d_{\cY}, \nu)$, the {\em Gromov-Wasserstein} (GW) distance between $\mathbb{X}$ and $\mathbb{Y}$ is defined by \cite{Mem:11}:
\begin{equation} \label{eq:gw_defn}
\gw_{p,q} (\mathbb{X},\mathbb{Y}) := \underset{\pi \in \Pi (\mu,\nu)}{\inf} \left( \int_{\cX \times \cY} \int_{\cX \times \cY} | d_{\cX} (\bx_1,\bx_2)^q - d_{\cY} (\by_1,\by_2)^q |^p ~ d \pi (\bx_1,\by_1) d \pi (\bx_2,\by_2) \right)^{1/p}.
\end{equation}
Here, the constraint set $\Pi (\mu,\nu)$ represents all feasible transportation plans -- these are modelled as joint probability measures over the space $\cX \times \cY$ subject to the condition that these satisfy the marginal constraints:
\begin{equation} \label{eq:marginal}
\int_{\cY} d \pi (\bx,\by) = d \mu, \quad \int_{\cX} d \pi (\bx,\by) = d \nu.
\end{equation}
At an intuitive level, the formulation in \eqref{eq:gw_defn} favours mappings $\pi$ that agree on the metric evaluations between pairs of points in $\cX$ with pairs of points in $\cY$.  The GW distance also confers certain geometric information -- if one is to consider the collection of all possible metric spaces, then the GW distance also defines a pseudo-metric over this collection, sometimes referred to as the $L^{p,q}$ {\em distortion distance} \cite{Sturm:23}.

\subsection{The Discrete GW problem and the Moment-SOS hierarchy}

From an optimization perspective, the formulation that is most convenient is the following:
\begin{equation} \label{eq:gw}
\gwp (\mathbb{X},\mathbb{Y}) := \underset{\pi \in \Pi (\mu,\nu)}{\inf} \int_{\cX \times \cY} \int_{\cX \times \cY} c(\bx_1,\by_1,\bx_2,\by_2)~ d \pi (\bx_1,\by_1) d \pi (\bx_2,\by_2).
\end{equation}
To be clear, the definition in \eqref{eq:gw_defn} is equivalent to making the following choice of cost function $c$ within \eqref{eq:gw}, and subsequently computing the $p$-th root: 
\begin{equation} \label{eq:lp_definition}
c(\bx_1,\by_1,\bx_2,\by_2) := | d_{\cX} (\bx_1,\bx_2)^q - d_{\cY} (\by_1,\by_2)^q |^p.
\end{equation}

In practice, one goes about solving \eqref{eq:gw} by discretizing the spaces $\CX$ and $\CY$, from which we obtain a finite dimensional optimization instance.  More concretely, suppose that $\CX = \{\bx_1,\dots,\bx_m\}$ and $\CY = \{\by_1,\dots,\by_n\}$ respectively.  Set $\mu(\bx_i) = \mu_i$ and $\nu(\by_j) =\nu_j$ for all $i \in [m]$, $j \in [n]$.  Note in particular that all feasible {\em couplings} (solutions) $\pi$ are now supported over the set $\{(\bx_i,\by_j)\;:\; i \in [m],\; j\in [n]\}$.  As such, we denote $\pi_{ij} : = \pi( (\bx_i,\by_j))$ for all $i \in [m]$, $j\in [n]$.  Then the collection of all feasible coupling measures is specified by the following polyhedral set:
\begin{equation*}
\begin{aligned}
\Pi(\mu,\nu) = \bigl\{ \pi \in \RR^{m \times n}~:~ & e_{ij}(\pi) := \pi_{ij} \geq 0 ~ \forall i \in [m],\; j \in [n], \\
& r_i(\pi) := \textstyle \sum_{j \in [n]} \pi_{ij} - \mu_i =0 ~ \forall i \in [m], \\
& c_j(\pi):= \textstyle \sum_{i \in [m]} \pi_{ij} - \nu_j =0 ~ \forall j \in [n] \bigr\}.
\end{aligned}
\end{equation*}
In particular, the discrete GW problem can be expressed as the following quadratic programming instance:
\begin{equation} \label{eq:gw_pop}
\min ~~ \sum_{i,j,k,l} c_{ij,kl}\cdot \pi_{ij}\pi_{kl} \qquad \mathrm{s.t.} \qquad \pi \in \Pi(\mu,\nu).
\end{equation}
Here, we denote $c_{ij,kl} = c(\bx_i,\by_j,\bx_k,\by_l)$ for all $i,k \in [m]$, and $j,l \in [n]$. 

Unfortunately, the optimization instance \eqref{eq:gw_pop} is challenging to solve.  First, the optimization instance in \eqref{eq:gw_pop} is a non-convex quadratic program.  In fact, recent work in \cite{kravtsova2024np} shows that the optimization formulation \eqref{eq:gw_pop} is NP-Hard.  The GW problem \eqref{eq:gw} is also conceptually connected to graph matching problems as well as the Quadratic Assignment Problem (QAP), both of which are known to be NP-Hard \cite{SG:76}.  In short, we simply do not expect to be able to solve \eqref{eq:gw_pop} using tractable algorithms.  

{\bf The moment-SOS hierarchy for the Discrete GW problem.} To address these difficulties, the authors in \cite{CNKS:24,tran:25} propose a series of tractable relaxations for solving \eqref{eq:gw_pop} using the ideas of the sum-of-squares (SOS) hierarchy.  Broadly speaking, the SOS hierarchy provides a principled framework for describing tractable relaxations to general polynomial optimization problems (POPs) \cite{thesisparrilo,lasserre2001global,BleParTho:12}.  It builds on two key ideas -- first, that SOS polynomials are globally non-negative, and therefore can be deployed as certificates (proofs) of global non-negativity of a given polynomial, and second, the fact that the problem of deciding if a polynomial is a SOS can be expressed in terms of a semidefinite program, which can be solved tractably \cite{Ren:01}.  To this end, the work in \cite{tran:25} describes a hierarchy of semidefinite relaxations to the discrete GW problem based on the moment-SOS hierarchy \cite{lasserre2001global, lasserre2009moments, lasserre2011new}.  The work proves several useful properties concerning these relaxations:  First, it shows that the optimal value of these relaxations converges to the global solution to the original discrete GW instance.  Second, it shows that these relaxations induce a type of pseudo-metric over the space of metric measure spaces, in a similar fashion that the GW distance does.

\subsection{Our contributions}

The moment-SOS hierarchy proposed in \cite{tran:25} is a hierarchy of semidefinite relaxations that aims to solve \eqref{eq:gw_pop}, which in reality is a discretization of the continuous problem in \eqref{eq:gw}.  A natural question is the following: is there a continuous extension of the ideas in \cite{tran:25}?  Can we view the proposed moment-SOS hierarchy as the discretized instance of a suitably defined continuous optimization instance that relaxes \eqref{eq:gw}?  What is the relationship between these problems?  Is there a physical interpretation behind the continuous analog?    

These are the main contributions of this paper:  In Section \ref{sec: moment for cts GW} we identify a suitable continuous limit of the moment-SOS hierarchy.   In Section \ref{sec:metric}, we show that the limit induces a pseudo-metric over the space of metric measure spaces that differs from the $L_{p,q}$ distortion distance.  In particular, the result relies on an interesting variant of the so-called gluing lemma.  The key insight we provide is that there is a suitable construction of the joint probability measure (the ``glued-measure'') that behaves well even under semidefinite constraints.  In Section \ref{sec:consistency} we prove a statistical consistency type of result in which data-points are drawn from the spaces $\CX$ and $\CY$ uniformly at random.  We prove that the optimal value of the discrete problem converges to the optimal value of the continuous problem as the number of samples increases.

Our intention in writing this paper is not merely about identifying the continuous limit of a moment-SOS hierarchy, so to speak; in fact, we have more ambitious goals.  The SOS hierarchy is traditionally applied to polynomial optimization instances.  An interesting question is if these techniques also apply to optimization problems over distributions.  The ideas outlined in this paper provide the agenda for the GW problem \eqref{eq:gw}, and key to this is the observation that the objective in \eqref{eq:gw} has a quadratic dependency on $\pi$.  Do our ideas also suggest a natural extension to general optimization problems that depend, in a polynomial way, on probability distributions?  We think this is an exciting question and we hope our work lays foundational ideas in this space.

\subsection{Notation and Outline}  We use $x$ to denote a scalar, and boldface $\bx$ to denote a vector.  Given $\bx = (x_1,\dots,x_n)$, we let $\RR[\bx]$ denote the polynomial ring over the variables $x_1,\ldots,x_n$.  We denote a monomial over the variables $\bx$ by $\bx^{\balpha} = \prod_{i=1}^nx_i^{\alpha_i}$, where $\balpha = (\alpha_1,\dots,\alpha_n) \in \NN^n$ is a multi-index of dimension $n$.  We denote the length of multi-index $\alpha$ by $|\balpha| = \sum_{i=1}^n \alpha_i$.  We denote the set of multi-indices with length at most $d$ by $\NN^n_d$, and we denote the set of multi-indices with length equal to $d$ by $\overline{\NN}^n_d$.  Given a positive integer $d$, we let $\RR[\bx]_d$ denote the set of polynomials of degree at most $d$.  The standard monomial basis is $\bv_d := (\bx^{\balpha})_{\balpha \in \NN^n_d}$, and has size $s(n,d) := \binom{n+d}{d}$.   We express polynomials $f \in \RR[\bx]$ as $f(\bx)= \sum_{\balpha \in \NN^n}f_{\balpha}\bx^{\balpha}$.  Let $\cK \subset \RR^n$.  We denote the collection of probability measures over $\cX$ by $\cP(\cX)$.  Last, we let $\biy = (\iy_{\balpha})_{\balpha \in \NN^n}$ denote moment sequences -- these are sequences of reals indexed by monomials $\balpha$.

In Section \ref{sec:momentsos} we describe the moment-SOS hierarchy for \eqref{eq:gw} based on the ideas in \cite{tran:25}.  In Section \ref{sec:mommeasure}, we describe an alternative formulation of the moment-SOS hierarchy.  The latter formulation lends itself to generalizations to continuous measures $\mu$ and $\nu$, which we describe in Section \ref{sec: moment for cts GW}.  In Section \ref{sec:metric} we discuss the metric properties arising from the distance measures defined in Section \ref{sec: moment for cts GW}, and in Section \ref{sec:consistency} we discuss statistical consistency questions arising from sampling $\mu$ and $\nu$ randomly.

To avoid confusion, we remind the reader that the optimization formulations \eqref{eq:gw_defn} and \eqref{eq:gw} are equivalent provided one chooses $c$ as in \eqref{eq:lp_definition}.  In Sections \ref{sec:momentsos}, \ref{sec:mommeasure}, and \ref{sec: moment for cts GW}, we appeal to the formulation \eqref{eq:gw} to emphasize that these ideas are relevant for {\em general} cost functions, and not simply those of the form \eqref{eq:lp_definition}.  Later on, in Sections \ref{sec:metric} and \ref{sec:consistency}, we appeal to the formulation \eqref{eq:gw_defn} because we need certain topological structure that arises from the specific choice of cost function \eqref{eq:lp_definition}.

\section{The Moment-SOS hierarchy for the Discrete GW} \label{sec:momentsos}

In this section, we describe the moment-SOS relaxation for the discrete GW problem \eqref{eq:gw_pop} proposed in \cite{tran:25}.  We begin our description by providing a brief background concerning the more general problem of minimizing a polynomial $f$ over a semialgebraic set $\cK$:
\begin{equation}\label{eq:pop}
\fmin := \inf \{ f(\bx) : \bx \in \cK \} \quad \text{where} \quad \cK = \{\bx \in \RR^n : g_j(\bx) \geq 0 \; \forall j \in [k] \}.
\end{equation}
Here, the $g_j$'s are polynomials in $\bx$.  The POP instance \eqref{eq:pop} can be re-formulated as a different optimization instance in which the variables are moment valuations.  Let $\biy = (\iy_{\balpha})_{\balpha \in \NN^n}$ be a sequence of reals indexed by the multi-indices $\balpha$.  Then \eqref{eq:pop} is equivalent to the following:
\begin{equation}\label{eq:pop2}
\min_{y_{\balpha}} ~~ \sum_{\balpha \in \NN^n}f_{\balpha} y_{\balpha} \quad \mathrm{s.t.} \quad y_{\balpha} = \mathbb{E}_{\bP} [\bx^{\balpha}] \text{ for some probability measure } \bP .
\end{equation}
Stated in words, \eqref{eq:pop2} minimizes over the values of each entry of the moment sequence $\biy$, subject to the constraint that these values are indeed realizable by some probability measure over $\cK$.  At optimality, the probability measure corresponding to the optimal moment sequence $\biy^{\star}$ is (almost surely) supported over the optimal values of $\bx$ in \eqref{eq:pop}.

Our next task is to describe the set of all feasible moment sequences $\biy$ derived from some measure $\bP$.  First, to each $\biy$, we can associate to it the {\em Riesz functional} $\ell_{\biy} : \RR[\bx] \to \RR$ such that the $\balpha$-component of $\biy$ satisfies: 
\begin{equation*}
y_{\balpha} = \ell_{\biy}(\bx^{\balpha})
\end{equation*}
so that
\begin{equation*}
\ell_{\biy}(f) = \sum_{\balpha \in \NN^n}f_{\balpha}\iy_{\balpha} \qquad \text{ for all } \qquad f(\bx) := \sum_{\balpha \in \NN^n}f_{\balpha}\bx^{\balpha}.
\end{equation*}

Given a positive integer $r \in \NN$, we define the {\em pseudo-moment matrix} of order $r$ associated with $\biy$ by: 
\begin{equation*}
\bM_r(\biy)(\balpha,\bbeta) := \ell_{\biy}(\bx^{\balpha+\bbeta})= \iy_{\balpha+\bbeta}, \quad \balpha,\bbeta \in \NN^n_r.
\end{equation*}
Similarly, given a polynomial $g \in \RR[\bx]$, we define the {\em localizing matrix} $\bM_r(g\biy)$ of order $r$ associated with $g$ and $\biy$ as the following:
\begin{equation*}
\bM_r(g\biy)(\balpha,\bbeta) := \ell_{\biy}(g(\bx)\bx^{\balpha+\bbeta})= \sum_{\bgamma \in \NN^n }g_{\bgamma}\iy_{\balpha+\bbeta+\bgamma}, \quad \forall \balpha, \bbeta \in \NN^n_{r- \lceil \deg (g) /2 \rceil}.
\end{equation*}
The {\em Riesz-Haviland} \cite[Theorem 3.1]{lasserre2009moments} states that a necessary and sufficient condition for $\biy$ to be realizable as the moment sequence corresponding to some probability measure is that the moment matrices $\bM_r(\biy)$ as well as the localizing matrices $\bM_r(g_i\biy)$ are PSD for all integers $r \in \NN$.  Unfortunately, this is a description that involves an {\em infinite} number of linear matrix inequalities.  To this end, the moment-SOS hierarchy is a relaxation that is derived by only requiring that the associated moment matrices and localizing matrices of order up to some level $r \leq r'$ be PSD, as opposed to all orders.  In particular, the latter relaxed description can be efficiently specified via semidefinite programs (SDPs).

More formally, the set of SOS polynomials is defined by:
\begin{equation*}
\Sigma [\bx] ~:=~ \Big \{ \sum_{i=1}^k p_i^2 : p_i \in \RR[\bx] \Big \}.
\end{equation*}
In a similar fashion, we let $ \Sigma[\bx]_{2r} \subset \Sigma [\bx]$ denote the subset of SOS polynomials with degree at most $2r$.  

The {\em pre-ordering} $\cT(\cK)$ is a subset of $\mathbb{R}[\bx]$ defined as follows:
\begin{equation*}
\cT(\cK) ~ :=~ \Big\{\sum_{j=1}^N\sigma_{I_j}(x)g_{I_j}(x) \in \RR[\bx] \,:\, \exists N \in \NN,\ \sigma_{I_j} \in \Sigma[x],\ g_I = \prod_{i \in I}g_i,\ g_{\emptyset}=1 \Big\}.
\end{equation*}
Similarly, the corresponding degree-bounded analog of the pre-ordering is given by:
\begin{equation*}
\cT(\cK)_{2r} ~:=~ \Big\{\sum_{j=1}^N\sigma_{I_j}(x)g_{I_j}(x) \in \cT(\cK) \,:\, \deg \sigma_{I_j}g_{I_j} \leq 2r \ \forall j \in [N] \Big\}.
\end{equation*}
Here, we denote $d_j := \lceil \deg g_j /2 \rceil$.  

Given an integer $r \in \NN$, the $r$-th level of the {\em Schmüdgen-type} moment hierarchy is given by:
\begin{equation} \label{eq:general_sch}
\begin{aligned}
\mlb(f,\mathcal{T}(\cK))_r ~:=~ \underset{\biy \in \RR^{s(n,2r)}}{\min} \quad & \sum_{\balpha \in \NN_{r}^n} f_{\balpha}\iy_{\balpha}\\
\mathrm{s.t.} \quad & \iy_0=1,\; \bM_{r-d_I}(g_I\biy) \succeq 0 \, \forall \  \{ I : I \subset [k], \deg g_I \leq 2r \}
\end{aligned}.
\end{equation}
Here, we denote $d_I := \lceil \deg g_I /2 \rceil$.  Note that the (constant) polynomial $1$ has degree $0$ and hence the above constraints automatically include the constraint $\bM_r(\biy) \succeq 0$.  The constraint $\bM_{r-d_I}(g_I\biy) \succeq 0$ is a finite dimensional linear matrix inequality.  In particular, \eqref{eq:general_sch} specifies an instance of a semidefinite program (SDP).  A consequence of the {\em Schmüdgen Positivstellensatz} is that, in the event where the set $\cK$ is compact, the sequence of lower bounds $\{\mlb(f,\mathcal{T}(\cK))_r\}_{r}^{\infty}$ increases monotonically to $f_{\min}$ as $r \rightarrow \infty$ \cite{schmudgen2017moment}.  

{\bf The Moment-SOS hierarchy for the Discrete GW.}  Based on these ideas, the work in \cite{tran:25} proposes the following Schmüdgen-type moment-SOS hierarchy for the discrete GW problem in \eqref{eq:gw_pop}:
\begin{equation} \label{eq:schmudgen_dgw}
\begin{aligned}
\underset{\biy \in \RR^{s(mn,2r)}}{\min} \qquad & \sum_{i,k \in [m],\; j,l \in [n]} c_{ij,kl}\cdot \ell_{\biy}(\pi_{ij}\pi_{kl}) \\
\mbox{subject to} \qquad & \iy_0 =1, \;\\ 
& \bM_{r - d_I}(e_I\biy) \succeq 0 \; \forall e_I = \textstyle \prod_{(i,j) \in I} e_{ij}, I \subset [m \times n],\; \deg e_I \leq 2r,\\
& \bM_{r-1}(r_i(\pi) \biy) = 0 \; \forall i \in [m],\\
& \bM_{r-1}(c_j(\pi) \biy) = 0 \; \forall j \in [n].
\end{aligned} \tag{S-DGW-r}
\end{equation}
The description in \eqref{eq:schmudgen_dgw} is a slight modification of the prescribed hierarchy specified in \eqref{eq:general_sch}.  Specifically, there is the appearance of the marginal conditions $\bM_{r-1}(r_i(\pi) \biy) = 0$ and $\bM_{r-1}(c_j(\pi) \biy) = 0$.  These are valid equalities and have a probabilistic interpretation of marginalization.  The inclusion of these valid inequalities turn out to be very useful as they allow one to obtain a different but substantially simpler {\em Putinar-type} moment-SOS hierarchy.  Specifically, it is possible to show that it suffices to only include the constraints of the type $\bM_{r - d_I}(e_I\biy) \succeq 0$ whereby $e_I = \textstyle \prod_{(i,j) \in I} e_{ij}$ is a {\em single} term, as opposed to all possible products (with degree at most $2r$) as in \eqref{eq:schmudgen_dgw} \cite{tran:25}.  

For the purpose of this work, it is more convenient to work with a slight re-formulation of \eqref{eq:schmudgen_dgw}.  Specifically, given $I$, let $\overline{\bM}_{r - d_I}(e_I\biy)$ be the principle submatrix of $\bM_{r - d_I}(e_I\biy)$ obtained by removing all columns and rows indexed by monomials $\pi^{\bgamma}$ whereby $|\bgamma| < r-d_I$.  Then the formulation \eqref{eq:schmudgen_dgw} is equivalent to the following \cite{tran:25}:
\begin{equation} \label{eq:schmudgen_dgw2}
\begin{aligned}
\underset{\biy \in \RR^{s(mn,2r)}}{\min} \qquad & \sum_{i,k \in [m],\; j,l \in [n]} c_{ij,kl}\cdot \ell_{\biy}(\pi_{ij}\pi_{kl}) \\
\mbox{subject to} \qquad & \iy_0 =1, \;\\ 
& \overline{\bM}_{r - d_I}(e_I\biy) \succeq 0 \; \forall e_I = \textstyle \prod_{(i,j) \in I} e_{ij}, I \subset [m \times n],\; \deg e_I = 2d \leq 2r,\\
& \bM_{r-1}(r_i(\pi) \biy) = 0 \; \forall i \in [m],\\
& \bM_{r-1}(c_j(\pi) \biy) = 0 \; \forall j \in [n].
\end{aligned} \tag{S-DGW-r2}
\end{equation}
 Note that in the above description, we only specify using subsets $I$ with {\em even} cardinality.  We explain why this is sufficient.  Suppose $I = \{ (i_1,j_1),\dots,(i_{2d-1},j_{2d-1})\}$ has odd cardinality.  The presence of the marginal conditions $\bM_{r-1}(r_i(\pi) \biy) = 0$ and $\bM_{r-1}(c_j(\pi) \biy) = 0$ imply
\begin{equation*}
\overline{\bM}_{r-d}(e_I\biy) = \textstyle \sum_{i \in [m],\; j \in [n]}\overline{\bM}_{r-d}((e_I\pi_{ij})\biy).
\end{equation*}
If $(i,j) \notin I$, then $I \cup {(i,j)}$ has even cardinality, and 
the term $\overline{\bM}_{r-d}((e_I\pi_{ij})\biy)$ is PSD because it appears among the constraints $\overline{\bM}_{r - d_I}(e_I\biy) \succeq 0$.  If instead $(i,j) \in I$, then $e_I \pi_{ij} = e_{I \backslash (i,j)} \pi_{ij}^2 $.  Then $\overline{\bM}_{r-d}((e_I\pi_{ij})\biy) = \overline{\bM}_{r-d}((e_{I \backslash (i,j)} \pi_{ij}^2)\biy)$, which appears as a principal submatrix of $\overline{\bM}_{r-d}((e_{I \backslash (i,j)})\biy)$ where $e_{I \backslash (i,j)}$ has even cardinality, and hence is also enforced to be PSD.

From a conceptual perspective, the formulation \eqref{eq:schmudgen_dgw} is simpler to understand; however, for analytical reasons, we refer to the equivalent formulation \eqref{eq:schmudgen_dgw2} as the definitive version of the moment-SOS hierarchy for the discrete GW problem \eqref{eq:gw_pop} in our subsequent results.  Finally, as a brief note, the work in \cite{tran:25} discusses an equivalent Putinar-type scheme (as mentioned earlier).  For the purposes of this paper, the ideas of Putinar-type scheme are not needed, and we shall make no further mention of this alternative scheme.

\section{From Moments to Measures} \label{sec:mommeasure}

The SOS hierarchy provided in \eqref{eq:schmudgen_dgw} is based on moment sequences.  Unfortunately, it is not immediately clear how the prescription generalizes to settings where the distributions $\mu,\nu$ may be continuous.  In what follows, we provide an alternative but equivalent description of \eqref{eq:schmudgen_dgw} whereby the decision variable is a {\em probability measure} rather than a moment sequence, and on which we may easily identify the continuous analog.  We start by describing the re-formulation over finite metric spaces, before progressing to the continuous extensions in the next section.  

In this section, we take $\CX = \{\bx_1,\dots,\bx_m\}$ and $\CY = \{\by_1,\dots,\by_n\}$.  Let $r \geq 1$ be a positive integer.  The optimization variable are probability measures over $(\CX \times \CY)^{2r}$; that is
\begin{equation*}
\bP \in \cP ((\CX \times \CY)^{2r}).
\end{equation*}
Because $\CX$ and $\CY$ are finitely supported, it means that $\bP \in (\mathbb{R}^{m \times n})^{\otimes 2r}$ is rank $2r$ tensor.  We refer to the $(i_1,j_1,\ldots,i_{2r},j_{2r})$-th coordinate of $\bP$ as follows:
\begin{equation*}
\bP ((\bx_{i_1},\by_{j_1},\dots,\bx_{i_{2r}},\by_{j_{2r}})) \in \mathbb{R} \quad \text{where} \quad i_s \in [m],\; j_s \in [n].
\end{equation*}

Let $S_n$ denote the collection of all permutations over $\{1,\ldots,n\}$.  The first constraint we impose is that $\bP$ is symmetric with respect to permutations over the indices:
\begin{equation}\label{eq:cond_sym}\tag{Sym}
\bP ((\bx_{i_1},\by_{j_1},\dots,\bx_{i_{2r}},\by_{j_{2r}}))= \bP ((\bx_{i_{\sigma(1)}},\by_{j_{\sigma(1)}},\dots,\bx_{i_{\sigma(2r)}},\by_{j_{\sigma(2r)}})) \quad \text{ for all } \quad \sigma \in S_{2r}.
\end{equation}

The second constraint we impose is that $\bP$ satisfies a type of marginal constraint
\begin{equation} \label{eq:cond_mar} \tag{Mar}
\bP |_{\cX \times (\cX \times \cY)^{2r-1}} = \mu \otimes \bP |_{(\cX \times \cY)^{2r-1}}, \quad \bP |_{\cY \times (\cX \times \cY)^{2r-1}}= \nu \otimes \bP |_{(\cX \times \cY)^{2r-1}}.
\end{equation}
As a reminder about notation, $\bP |_{\cX \times (\cX \times \cY)^{2r-1}}$ denotes the marginal distribution over $\cX \times (\cX \times \cY)^{2r-1}$ -- one way to understand this is to integrate over the indices that do not appear in the subscript.  For instance, in this example, we have the following:
\begin{equation*}
\bP |_{\bx_{i_1},\bx_{i_2},\by_{i_2},\dots,\bx_{i_{2r}},\by_{i_{2r}}} = \int_{ \by_{i_1} \in \CY } d\bP.
\end{equation*}
Because of \eqref{eq:cond_sym}, it is not critical which subset of indices we integrate over.

The third constraint we impose is that $\bP$ satisfies the following positive semidefinite condition:  Let $I = \{(i_1^*,j_1^*),\dots,(i_{2d_I}^*,j_{2d_I}^*) \} \subset [m \times n]$ be any subset of coordinates with $d_I \leq r$.  Let $f \in \cB((\cX \times \cY)^{r-d_I})$ be any measurable function.  We require the following to hold
\begin{equation}\label{eq:cond_psd} \tag{PSD}
\begin{aligned}
\int_{(\cX \times \cY)^{2r}}f(\overline{\bx}_1,\overline{\by}_1\dots,\overline{\bx}_{r-d_I},\overline{\by}_{r-d_I})f(\overline{\bx}_{r-d_I+1},\overline{\by}_{r-d_I+1},\dots,\overline{\bx}_{2r-2d_I},\overline{\by}_{2r-2d_I})\\
\cdot\prod_{s=1}^{2d_I} \delta_{(\overline{\bx}_{s+2r-2d_I},\overline{\by}_{s+2r-2d_I}) =(\bx_{i_s^*},\by_{j_s^*})}~d \bP \geq 0.
\end{aligned}
\end{equation}
Here is an alternative and perhaps more intuitive definition:  The subset $I$ selects $2d_{I}$ coordinates.  As such, the remaining coordinates that were not chosen specifies a tensor of rank-$2(r-d_I)$.  Because of symmetry, we can identify the first set of $(r-d_I)$ coordinates with the second set of $(r-d_I)$ coordinates; in particular, by vectorization, we may view this as a symmetric matrix of dimension $(mn)^{r-d_I} \times (mn)^{r-d_I}$.  The condition \eqref{eq:cond_psd} asks that this resulting matrix be PSD, for all possible choices of $I$.

The alternative moment-SOS hierarchy is given by the following semidefinite program:
\begin{equation} \label{eq:measure_SOSv1}
\begin{aligned}
\gwp^{(r)}(\mu,\nu) ~~:=~~ \inf_{\bP} \quad & \int_{(\cX \times \cY)^{2r}}c(\bx_1,\by_1,\bx_2,\by_2)~d\bP \\
\mathrm{s.t.} \quad & \bP \in \cP(\cX \times \cY)^{2r}, \\
& \bP \text{ satisfies } \eqref{eq:cond_sym}, \eqref{eq:cond_mar}, \eqref{eq:cond_psd}
\end{aligned}
\end{equation}

\subsection{Correspondence between \eqref{eq:schmudgen_dgw2}  and \eqref{eq:measure_SOSv1}}

Next, we show that \eqref{eq:schmudgen_dgw2} and \eqref{eq:measure_SOSv1} are indeed equivalent.  To do so, we show that feasible solutions in one problem can be mapped to feasible solutions in the other, with equal objective value.  In what follows, we let $\cM^{(r)}(\mu,\nu)$ denote the collection of moment sequences $\biy$ that are feasible in \eqref{eq:schmudgen_dgw2}.

[Mapping from $\cM^{(r)}(\mu,\nu)$ to $\cP((\cX \times \cY)^{2r})$]:  Consider the following mapping $\Py :\cM^{(r)}(\mu,\nu) \rightarrow \cP((\cX \times \cY)^{2r})$ specified as follows:
\begin{equation*}
\Py(\biy) [(\bx_{i_1},\by_{j_1},\dots,\bx_{i_{2r}},\by_{j_{2r}})] := \ell_{\biy} \Big (\prod_{s=1}^{2r}\pi_{i_sj_s} \Big) \quad \forall i_s \in [m],\; j_s \in [n].
\end{equation*}
First, we explain why $\ell_{\biy}(\prod_{s=1}^{2r}\pi_{i_sj_s} ) \geq 0$.  We can write $\prod_{s=1}^{2r}\pi_{i_sj_s}= e_I\pi^{2\alpha}$ for some $I \subset [m \times n]$ and some $\alpha \in \NN^n_r$ ($2\alpha$ collects the even exponents, and $I$ results from the odd exponents).  Then $\ell_{\biy} (\prod_{s=1}^{2r}\pi_{i_sj_s})$ corresponds to the $(\alpha,\alpha)$-entry of $\bM_{r-d_I}(e_I\biy)$.  Since $\bM_{r-d_I}(e_I\biy) \succeq 0$, it follows that $\ell_{\biy} (\prod_{s=1}^{2r}\pi_{i_sj_s}) \geq 0$.  Later, we verify that a marginalization-type condition ensures that the sum of these measures equals one, from which we conclude $\pi$ indeed specifies a probability measure.  Second, we briefly mention that there are no issues with well-definedness.  To see why, notice that $\CX$ and $\CY$ are finitely supported.  Hence, up to set of measure zero, the space is a finite collection of points.  Given any Borel subset $B \subset (\cX \times \cY)^{2r}$, we break it into a finite collection of atoms of the form $(\bx_{i_1},\by_{j_1},\dots,\bx_{i_{2r}},\by_{j_{2r}})$, from which we read off its probability mass via $\pi(\biy)$.

Next, we proceed with the constraints.  First, it is clear from construction that $\pi(\biy)$ satisfies \eqref{eq:cond_sym} for all $\biy \in \cM^{(r)}(\mu,\nu)$, and we omit the proof of this step.  Next:

\begin{lemma}\label{lem:mar}
Suppose $\biy \in \cM^{(r)}(\mu,\nu)$.  Then $\Py_{\biy}$ satisfies \eqref{eq:cond_mar}.
\end{lemma}

\begin{proof}[Proof of Lemma \ref{lem:mar}]
Because of the symmetry of $\pi_{\biy}$, for any $0 \leq k,l,t \leq 2r$ such that $k+l+t \leq 2r$ we have
\begin{equation*}
\begin{aligned}
& \Py_{\biy}((\bx_{i_1},\dots,\bx_{i_k},\by_{j_1},\dots,\by_{j_l},\bx_{i'_1},\by_{j'_1},\dots,\bx_{i'_t},\by_{j'_t})) \\ 
= ~ & \Py_{\biy}((\bx_{i_1}\times \cY,\dots,\bx_{i_k}\times \cY,\cX \times \by_{j_1},\dots,\cX \times \by_{j_l},\bx_{i'_1},\by_{j'_1},\dots,\bx_{i'_t},\by_{j'_t}, (\cX \times \cY)^{2r-t-k-l})).
\end{aligned}
\end{equation*}
Then, by the marginal conditions of $\biy$, for all $ 0 \leq t \leq 2r$, one has
\begin{align*}
&\Py_{\biy}(\bx_{i_1},\by_{j_1},\dots,\bx_{i_t},\by_{j_t})= \Py_{\biy}((\bx_{i_1},\by_{j_1},\dots,\bx_{i_t},\by_{j_t}) \times (\cX \times \cY)^{2r-t})\\
=\;& \sum_{\substack{i_{t+1},\dots,i_{2r} \in [m] \\j_{t+1},\dots,j_{2r} \in [n] }}\pi_{\biy}(\bx_{i_1},\by_{j_1},\dots,\bx_{i_{2r}},\by_{j_{2r}}) ~=~  \sum_{\substack{i_{t+1},\dots,i_{2r} \in [m] \\j_{t+1},\dots,j_{2r} \in [n] }}\ell_{\biy}\Big(\prod_{s=1}^{2r}\pi_{i_sj_s} \Big) =  \ell_{\biy}\Big( \prod_{s=1}^{t}\pi_{i_sj_s} \Big).
\end{align*}
Then we observe that
\begin{align*}
& \Py_{\biy}((\bx_{i_1},\bx_{i_2},\by_{j_2},\dots,\bx_{i_{2r}},\by_{j_{2r}})) = \Py_{\biy}((\bx_{i_1}\times \cY,\bx_{i_2},\by_{j_2},\dots,\bx_{i_{2r}},\by_{j_{2r}}))\\
=\;& \sum_{j_1 \in n} \Py_{\biy}((\bx_{i_1},\by_{j_1},\dots,\bx_{i_{2r}},\by_{j_{2r}}))= \sum_{j_1 \in n} \ell_{\biy} \Big(\prod_{s=1}^{2r}\pi_{i_sj_s} \Big)\\
=\;& \mu_{i_1}\ell_{\biy} \Big(\prod_{s=1}^{2r}\pi_{i_sj_s} \Big)= \mu_{i_1}\pi((\bx_{i_2},\by_{j_2},\dots,\bx_{i_{2r}},\by_{j_{2r}})).
\end{align*}
This implies the marginal condition 
\begin{equation*}
\Py_{\biy}|_{\cX \times (\cX \times \cY)^{2r-1}}= \mu \otimes \Py_{\biy}|_{(\cX \times \cY)^{2r-1}}.
\end{equation*}
A similar sequence of steps proves the other marginal condition.
\end{proof}

\begin{lemma}\label{lem: PSD}
Suppose $\biy \in \cM^{(r)}(\mu,\nu)$.  Then $\Py_{\biy}$ satisfies \eqref{eq:cond_psd}.
\end{lemma}

\begin{proof}[Proof of Lemma \ref{lem: PSD}]  The essence of this proof is to realize that \eqref{eq:cond_psd} exactly encodes the condition $\overline{\bM}_{r - d_I}(e_I\biy) \succeq 0$.  

Recall that $\Py_{\biy}$ is an atomic measure over $(\cX \times \cY)^{2r}$.  Our first step is to re-write the integral in \eqref{eq:cond_psd} as a sum.  Fix a subset $I \subset [m \times n]$, $I=\{ (i_1^*,j_1^*),\dots,(i_{2d_I}^*,j_{2d_I}^*) \}$.  Then, for any measurable function $f$ on $(\cX \times \cY)^{r}$, one has
\begin{equation} \label{eq:lemma_psd_intermediatestep}
\begin{aligned}
&\int_{(\cX \times \cY)^{2r}}f(\overline{\bx}_1,\overline{\by}_1\dots,\overline{\bx}_{r-d_I},\overline{\by}_{r-d_I})f(\overline{\bx}_{r-d_I+1},\overline{\by}_{r-d_I+1},\dots,\overline{\bx}_{2r-2d_I},\overline{\by}_{2r-2d_I})\\
& \qquad\qquad \cdot \prod_{s=1}^{2d_I}\delta_{(\overline{\bx}_{2r-2d_I+s},\overline{\by}_{2r-2d_I+s}) =(\bx_{i_s^*},\by_{j_s^*})}~d \Py_{\biy}\\
=\;& \sum_{\substack{i_1,\dots,i_{2r-2d_I} \in [m]\\j_1,\dots,j_{2r-2d_I} \in [n] }}f(\overline{\bx}_{i_1},\overline{\by}_{j_1},\dots,\overline{\bx}_{i_{r-d_I}},\overline{\by}_{j_{r-d_I}})f(\overline{\bx}_{i_{r-d_I+1}},\overline{\by}_{j_{r-d_I+1}},\dots,\overline{\bx}_{i_{2r-2d_I}},\overline{\by}_{j_{2r-2d_I}})\\
& \qquad\qquad \cdot \Py_{\biy}((\overline{\bx}_{i_1},\overline{\by}_{j_1},\dots,\overline{\bx}_{i_{2r-2d_I}},\overline{\by}_{j_{2r-2d_I}},\overline{\bx}_{i^*_1},\overline{\bx}_{j^*_1},\dots,\overline{\bx}_{i^*_{2d_I}},\overline{\by}_{i^*_{2d_I}}))\\
=\;& \sum_{\substack{i_1,\dots,i_{2r-2d_I} \in [m]\\j_1,\dots,j_{2r-2d_I} \in [n] }}f(\overline{\bx}_{i_1},\overline{\by}_{j_1},\dots,\overline{\bx}_{i_{r-d_I}},\overline{\by}_{j_{r-d_I}})f(\overline{\bx}_{i_{r-d_I+1}},\overline{\by}_{j_{r-d_I+1}},\dots,\overline{\bx}_{i_{2r-2d_I}},\overline{\by}_{j_{2r-2d_I}})\\
& \qquad\qquad \cdot \ell_{\biy}\left(\prod_{s=1}^{2r-2d_I}\pi_{i_sj_s} \prod_{s=1}^{2d_I}\pi_{i^*_sj^*_s}\right).
\end{aligned}
\end{equation}

The remainder of this proof is to recognize that this expression is of the form $\mathbf{f}^{\top}\overline{\bM}_{r-d_I}(e_I\biy) \mathbf{f}$ of some vector $\mathbf{f}$.  Given a multi-index $\bgamma \in \NN^{mn}_{r-d_I}$ of length $r-d_I=|\bgamma| = \sum_{i \in [m],\; j\in [n]}|\gamma_{ij}|$, we define the function 
\begin{equation*}
f(\pi^{\bgamma}):= \sum_{\substack{i_1,\dots,i_{r-d_I} \in [m],\; j_1,\dots,j_{r-d_I} \in [n]\\ \prod_{s=1}^r\pi_{i_sj_s}=\pi^{\bgamma}}}f(\bx_{i_1},\by_{j_1},\dots,\bx_{i_{r-d_I}},\by_{j_{r-d_I}}).
\end{equation*}
We denote $\mathbf{f} = ( f(\pi^{\bgamma}) )_{\bgamma \in \NN^{mn}_{r-d_I,\; |\bgamma| = r-d_I}}$ whose coordinates are indexed by $\bgamma$, and whose $\bgamma$-th corodinate is $f(\pi^{\bgamma})$.  This allows us to re-write the sum in \eqref{eq:lemma_psd_intermediatestep} as 
\begin{equation*}
\sum_{\bgamma, \btau \in \NN^{mn}_{r-d_I}, |\bgamma| = |\btau| =r-d_I }f(\pi^{\bgamma})f(\pi^{\btau})\ell_{\biy}(\pi^{\bgamma+\btau}) ~=~ \mathbf{f}^{\top}\overline{\bM}_{r-d_I}(e_I\biy) \mathbf{f} ~\geq~ 0.
\end{equation*}
Here, we obtain the last inequality because $\overline{\bM}_{r-d_I}(e_I\biy) \succeq 0$.
\end{proof}

First, we note that the objective can be expressed as the following integral
\begin{equation*}
\sum_{i,k \in [m],\; j,l \in [n]} c_{ij,kl}\cdot \ell_{\biy}(\pi_{ij}\pi_{kl})= \int_{(\cX \times \cY)^{2r}}c(\overline{\bx}_1,\overline{\by}_1,\overline{\bx}_2,\overline{\by}_2)~d\pi_{\biy}.
\end{equation*}

[Mapping from $\cP((\cX \times \cY)^{2r})$ to $\cM^{(r)}(\mu,\nu)$]:  Next, we specify the reverse map.  Let $\bP \in \cP((\cX \times \cY)^{2r}$ be a probability distribution.  We begin by specifying the moments across all monomials of degree $2r$:
\begin{equation*}
\ell \Big (\prod_{s=1}^{2r}\pi_{i_s j_s} \Big) := \bP [(\bx_{i_1},\by_{j_1},\dots,\bx_{i_{2r}},\by_{j_{2r}})] \quad \forall i_s \in [m],\; j_s \in [n].
\end{equation*}
For all remaining moments $\pi^{\gamma}$ where $|\gamma| < 2r$, we specify the moments by
\begin{equation} \label{eq:reversemap_smallermoments}
\ell ( \pi^{\gamma} ) := \sum_{\tau : |\tau + \gamma| = 2r} \ell ( \pi^{\tau} \pi^{\gamma} ) .
\end{equation}

We need to explain why this construction leads to a moment sequence $\biy$ that is feasible in \eqref{eq:schmudgen_dgw}.  First, we have $\iy_0 =1$ because $\bP$ is a probability distribution.  Second, the marginal constraints $\bM_{r-1}(r_i(\pi) \biy) = 0$ and $\bM_{r-1}(c_j(\pi) \biy) = 0$ hold because of \eqref{eq:reversemap_smallermoments}.  Third, to see that the PSD constraint $\bM_{r - d_I}(e_I\biy) \succeq 0$ holds, we reverse the steps of Lemma \ref{lem: PSD}.

[Equal objective value]:  Last, we check that the objective values of \eqref{eq:schmudgen_dgw2} and \eqref{eq:measure_SOSv1} are equal.  This step is trivial -- because $\CX$ and $\CY$ are finite metric spaces, the integral in \eqref{eq:measure_SOSv1} converts to sums, and one can check they take the form in \eqref{eq:schmudgen_dgw2}.





\section{Moment-SOS Hierarchy for the Continuous GW}\label{sec: moment for cts GW}

The objective of this section is to state the extension of \eqref{eq:measure_SOSv1} so that it is compatible with continuous distributions $\mu$ and $\nu$.  As before, the main optimization variable is a probability measure $\bP \in \cP((\cX \times \cY)^{2r})$.  Because the distributions $\mu$ and $\nu$ are continuous, it is necessary to put in place additional assumptions, which we state as follows:

\begin{enumerate}
\item We assume that the metric spaces $(\CX,d_{\CX})$ and $(\CY,d_{\CY})$ are Polish and {\em compact}.  Compactness allows us to access finite covers later on.  We also assume all measures are Borel.
\item In particular, we take $\CX$ and $\CY$ have diameter one. 
\item We assume that the cost function $c(d_{\CX} (x_0,x_1), d_{\CY} (y_0,y_1))$ is Lipschitz continuous on $(\cX \times \cY) \times (\cX \times \cY)$ with respect to the distance $d_{\cX}+d_{\cY}$.
\end{enumerate}



First, we note that the symmetric constraints \eqref{eq:cond_sym} and \eqref{eq:cond_mar}, as it is stated, already make sense for continuous measures.  As such, there is no need to further generalize these beyond the minor qualification that we now only require \eqref{eq:cond_sym} and \eqref{eq:cond_mar} to hold almost surely.

%

Next, let $\cB_+((\cX \times \cY)^{n})$ denote the set of non-negative measurable functions over $(\cX \times \cY)^{n}$.  We extend the condition \eqref{eq:cond_psd} as follows:  For all $t$ such that $0 \leq t \leq r$, all measurable $f \in \cB((\cX \times \cY)^{2t})$, and all non-negative measurable $g \in \cB_+((\cX \times \cY)^{2r-2t})$, we ask that the following integral be non-negative:
\begin{equation*}\label{cond: PSD for cts} \tag{PSD+}
\int_{(\cX \times \cY)^{2r}} f(\bx_1,\by_1,\dots,\bx_t,\by_{t}) f(\bx_{t+1},\by_{t+1},\dots,\bx_{2t},\by_{2t}) g(\bx_{2t+1},\by_{2t+1},\dots,\bx_{2r},\by_{2r})~d\bP \geq 0 .
\end{equation*}

The proposed relaxation of \eqref{eq:gw} is as follows:
\begin{equation} \label{hierarchy moment for cst}
\begin{aligned}
\gwp^{(r)}  (\mu,\nu) := \quad & \inf_{\bP \in \cP(\cX \times \cY)^{2r}}\int_{(\cX \times \cY)^{2r}}c(\overline{\bx}_1,\overline{\by}_1,\overline{\bx}_2,\overline{\by}_2)~d\bP \\
\mathrm{s.t.} \quad & \bP \text{ satisfies } \eqref{eq:cond_sym}, \eqref{eq:cond_mar}, \eqref{cond: PSD for cts}.
\end{aligned}
\end{equation}


Let $\gwp^{(r)}(\mu,\nu)$ denote the optimal value of \eqref{hierarchy moment for cst}.  Let $\gwp(\mu,\nu)$ denote the optimal value of \eqref{eq:gw}.  How do these values relate to each other?  This is the subject of the following results.  First:

\begin{proposition} \label{thm:moment_lowerbounds}
One has $\gwp^{(r)}(\mu,\nu) \leq \gwp(\mu,\nu)$ for all $r$.
\end{proposition}

\begin{proof}[Proof of Proposition \ref{thm:moment_lowerbounds}]
Fix a positive integer $r$.  To simplify notation, we denote $\bz := (\bx,\by)$, and $\CZ := (\CX \times \CY)$.  Recall from \eqref{eq:gw}
\begin{equation*}
\gwp(\mu,\nu) := \underset{\pi \in \Pi (\mu,\nu)}{\inf} \int_{\CZ} \int_{\CZ} c(\bz_1,\bz_2) ~ d \pi (\bz_1) d \pi (\bz_2).
\end{equation*}

To show that $\gwp^{(r)}(\mu,\nu) \leq \gwp(\mu,\nu)$, it suffices to show that, for any $\pi \in \Pi(\mu,\nu)$, the measure $\pi^{2r}:= \otimes^{2r} \pi$ is a feasible solution to \eqref{hierarchy moment for cst}.  It is clear from the definition that \eqref{eq:cond_sym} and \eqref{eq:cond_mar} hold.  Let $t \leq r$ be a positive integer.  Given $f \in \cB(\CZ^{2t})$ and $ g \in \cB_+(\CZ^{2r-2t})$, we have the following:
\begin{equation*}
\begin{aligned}
&\int_{\CZ^{2r}}f(\bz_1,\dots,\bz_t)f(\bz_{t+1},\dots,\bz_{2t}) g(\bz_{2t+1},\dots,\bz_{2r}) ~ d\pi^{2r}\\
=\;& \big (\int_{\CZ^{t}}f(\bz_1,\dots,\bz_t)~d\pi^t \big ) \big (\int_{\CZ^{t}}f(\bz_{t+1},\dots,\bz_{2t})~d\pi^t \big ) \cdot \big (\int_{\CZ^{2r-2t}} g(\bz_{2t+1},\dots,\bz_{2r})~d\pi^{2r-2t} \big )\\
=\;& \big (\int_{\CZ^{t}} f(\bz_1,\dots,\bz_t)~d\pi^t \big )^2 \big (\int_{\CZ^{2r-2t}} g(\bz_{2t+1},\dots,\bz_{2r})~d\pi^{2r-2t} \big ) \geq 0.
\end{aligned}
\end{equation*}
The last inequality is due to the fact that $g$ is a non-negative measurable function on $\CZ^{2r-2t}$.  As such, \eqref{cond: PSD for cts} also holds, and hence $\pi^{2r}$ is feasible in \eqref{hierarchy moment for cst}.
\end{proof}

In fact, we can strengthen the result to:

\begin{theorem} \label{thm: convergence}
One has $\gwp^{(r)}(\mu,\nu) \rightarrow \gwp(\mu,\nu)$ as $r \rightarrow \infty$.
\end{theorem}


Our proof of convergence of the sequence $\{ \gwp^{(r)}(\mu,\nu) \}_{r \in \NN}$ relies on suitably approximating the continuous problem with a GW discretized instance, for which we can apply prior results concerning convergence.

Consider the compact metric spaces $(\cX,d_{\cX})$ and $(\cY,d_{\cY})$.  Because these are compact, the spaces $\cX$ and $\cY$ are also {\em totally bounded}; that is to say, for any $\varepsilon > 0$, there exists partitions of $\CX$ and $\CY$ into finitely many subsets
\begin{equation}\label{partitions}
\CX = \bigsqcup_{i=1}^k \CX_i, \quad \mathrm{and} \quad \CY = \bigsqcup_{j=1}^l \CY_i
\end{equation}
such that for any $i \in [k], j \in [l]$, there exists $\bx_i \in \CX _i$ and $\by_j \in \CY_j$ satisfying 
\begin{equation*}
\CX_i \subset B(\bx_i,\varepsilon), \quad \mathrm{and} \quad \CY_j \subset B(\by_j,\varepsilon).
\end{equation*}
Subsequently, for any positive integer $s$, the space $(\cX \times \cY)^s$ admits the following partition 
\begin{equation*}
(\cX \times \cY)^s= \bigsqcup_{\substack{i_1,\dots,i_s \in [k]\\ j_1,\dots,j_s \in [l]}}\prod_{t=1}^s(\cX_{i_t} \times \cY_{j_t}), 
\end{equation*}
and in which case, we also have these inclusions
\begin{equation*}
\prod_{t=1}^s(\cX_{i_t} \times \cY_{j_t}) \subset B((\bx_{i_1},\by_{j_1},\dots,\bx_{i_s},\by_{j_s}),2s\varepsilon).    
\end{equation*}

Let $C$ be the Lipschitz constant of the function $c$ over the metric space $((\cX \times \cY)^2, d_{\cX}+d_{\cX}+d_{\cY}+d_{\cY})$.  Then, for any probability measure $\kappa$ on $(\cX \times \cY)^2$, we have the following
\begin{equation}\label{eq: lipszhitz inequality}
\Big|\int_{(\cX \times \cY)^{2}}c(\overline{\bx}_1,\overline{\by}_1,\overline{\bx}_2,\overline{\by}_2)~d\kappa -\sum_{\substack{i_1,i_2 \in [k]\\ j_1,j_2 \in [l]}}c(\bx_{i_1},\by_{j_1},\bx_{i_2},\by_{j_2})\kappa\left(\cX_{i_1}\times \cX_{i_2}\times \cY_{j_1}\times \cY_{j_2} \right) \Big| \leq 4C \varepsilon.
\end{equation}

Let $\pi \in \Pi(\mu,\nu)$ be a coupling measure.  We define the $\varepsilon$-{\em concentration} $\pi^{2r}_{c,\varepsilon}$ (with respect to the above partition) as follows: $\pi^{2r}_{c,\varepsilon}$ is an atomic measure whose atoms are
\begin{equation*}
\{(\bx_{i_1},\by_{j_1},\dots,\bx_{i_{2r}},\by_{j_{2r}})\;:\; i_1,\dots,i_{2r} \in [m],\; j_1,\dots,j_{2r} \in [n]\},
\end{equation*}
and whose corresponding probability mass is specified by
\begin{equation*}
\pi^{2r}_{c,\varepsilon}((\bx_{i_1},\by_{j_1},\dots,\bx_{i_{2r}},\by_{j_{2r}})) := \pi^{2r}\Big (\prod_{t=1}^{2r}(\cX_{i_t} \times \cY_{j_t}) \Big )=\prod_{s=1}^{2r}\pi(\cX_{i_t} \times \cY_{j_t}).
\end{equation*}
Similarly, we define the $\varepsilon$-concentration on the marginal condition $\mu$ and $\nu$ as 
\begin{equation*}
\mu_{c,\varepsilon}(\bx_i) = \mu(\cX_i)\; \forall i \in [k], \quad \mathrm{and} \quad \nu_{c,\varepsilon}(\by_i) = \mu(\cY_i)\; \forall j \in [l].
\end{equation*}

\begin{lemma}\label{thm:concentrationbound}
We have the following inequality: 
\begin{equation}\label{eq:concentrationinequality}
\gwp(\mu_{c,\varepsilon},\nu_{c,\varepsilon}) \leq \gwp(\mu,\nu) + 4C\varepsilon.
\end{equation}
\end{lemma}

\begin{proof}[Proof of Lemma \ref{thm:concentrationbound}]
First, given $\pi \in \Pi(\mu,\nu)$, we define the $\varepsilon$-concentration given by
\begin{equation*}
\pi_{c,\varepsilon} ((\bx_i,\by_j))= \pi(\cX_i \times \cY_j) \; \forall i \in [k],\; j\in [l].
\end{equation*}
Then, by combining \eqref{eq: lipszhitz inequality} with the fact that $\pi_{c,\varepsilon} \in \Pi(\mu_{c,\varepsilon},\nu_{c,\varepsilon})$, we have
\begin{equation*} 
\int_{(\cX \times \cY)^{2}}c(\overline{\bx}_1,\overline{\by}_1,\overline{\bx}_2,\overline{\by}_2)~d\pi_{c,\varepsilon} \leq \int_{(\cX \times \cY)^{2}}c(\overline{\bx}_1,\overline{\by}_1,\overline{\bx}_2,\overline{\by}_2)~d\pi + 4C\varepsilon .
\end{equation*}
By minimizing over $\pi \in \Pi(\mu,\nu)$ we have obtained \eqref{eq:concentrationinequality}.
\end{proof}

In what follows, we state and prove the analogous conclusion to $\gwp^{(r)}$.

\begin{lemma}\label{thm:concentrationhierarchy}
Given $\bP \in \Pi_r(\mu,\nu)$, we define its $\varepsilon$-concentration by:
\begin{equation*}
\bP_{c,\varepsilon}((\bx_{i_1},\by_{j_1},\dots,\bx_{i_{2r}},\by_{j_{2r}}))= P\left(\prod_{t=1}^{2r}(\cX_{i_t} \times \cY_{j_t}) \right), \; \forall \; i_1,\dots,i_{2r} \in [k],\; j_1,\dots,j_{2r} \in [l].
\end{equation*}
Then $\bP_{c,\varepsilon} \in \Pi_r(\mu_{c,\varepsilon},\nu_{c,\varepsilon})$.  In particular, this implies
\begin{equation}\label{ineq 3}
\gwp^{(r)}(\mu_{c,\varepsilon},\nu_{c,\varepsilon}) \leq \gwp^{(r)}(\mu,\nu) + 4C \varepsilon.
\end{equation}
\end{lemma}

\begin{proof}[Proof of Lemma \ref{thm:concentrationhierarchy}]
It is straightforward to check that the conditions \eqref{eq:cond_sym} and \eqref{eq:cond_mar} hold, and we omit these proofs.  We proceed to check \eqref{cond: PSD for cts}.  To this end, let $0 \leq t \leq r$, $f \in \cB((\cX \times \cY)^{2t})$, and $g \in \cB_+((\cX \times \cY)^{2r-2t})$.  Define the following step functions 
\begin{equation*}
\begin{aligned}
\overline{f}(\overline{\bx}_1,\overline{\by}_1,\dots,\overline{\bx}_t,\overline{\by}_t) := & \sum_{\substack{i_1,\dots,i_t \in [k]\\ j_1,\dots,j_t \in [l]}}f(\bx_{i_1},\by_{j_1},\dots,\bx_{i_t},\by_{j_t}) \delta_{\prod_{s=1}^t(\cX_{i_s}\times \cY_{j_s})},\\
\overline{g}(\overline{\bx}_1,\overline{\by}_1,\dots,\overline{\bx}_{2r-2t},\overline{\by}_{2r-2t}) := & \sum_{\substack{i_1,\dots,i_{2r-2t} \in [k]\\ j_1,\dots,j_{2r-2t} \in [l]}}g(\bx_{i_1},\by_{j_1},\dots,\bx_{i_{2r-2t}},\by_{j_{2r-2t}}) \delta_{\prod_{s=1}^{2r-2t}(\cX_{i_s}\times \cY_{j_s})}.
\end{aligned}
\end{equation*}
Since $P$ is supported in the set $\{(\bx_{i_1},\by_{j_1},\dots,\bx_{i_{2r}},\by_{j_{2r}}) : i_1,\dots,i_t \in [k], j_1,\dots,j_t \in [l] \}$, the definition of $\bP_{c,\varepsilon}$ implies that 
\begin{align*}
&\int_{(\cX \times \cY)^{2r}}f(\overline{\bx}_1,\overline{\by}_1,\dots,\overline{\bx}_t,\overline{\by}_t)f(\overline{\bx}_{t+1},\overline{\by}_{t+1},\dots,\overline{\bx}_{2t},\overline{\by}_{2t})g(\overline{\bx}_{2t+1},\overline{\by}_{2t+1},\dots,\overline{\bx}_{2r},\overline{\by}_{2r})~d\bP_{c,\varepsilon}\\
=\;& \int_{(\cX \times \cY)^{2r}}\overline{f}(\overline{\bx}_1,\overline{\by}_1,\dots,\overline{\bx}_t,\overline{\by}_t)\overline{f}(\overline{\bx}_{t+1},\overline{\by}_{t+1},\dots,\overline{\bx}_{2t},\overline{\by}_{2t})\overline{g}(\overline{\bx}_{2t+1},\overline{\by}_{2t+1},\dots,\overline{\bx}_{2r},\overline{\by}_{2r})~d\bP \geq 0.
\end{align*}
The last inequality comes from the positive semidefiniteness of $P$.  This means that \eqref{cond: PSD for cts} is satisfied. Hence, $\bP_{c,\varepsilon} \in \Pi_r(\mu_{c,\varepsilon},\nu_{c,\varepsilon})$.  By applying \eqref{eq: lipszhitz inequality}, we obtain
\begin{equation*}
\int_{(\cX \times \cY)^{2r}}c(\overline{\bx}_1,\overline{\by}_1,\overline{\bx}_2,\overline{\by}_2)~d\bP_{c,\varepsilon} \leq \int_{(\cX \times \cY)^{2r}}c(\overline{\bx}_1,\overline{\by}_1,\overline{\bx}_2,\overline{\by}_2)~d\bP + 4C\varepsilon.
\end{equation*}
By minimizing over $\bP \in \Pi_r(\mu,\nu)$, we conclude $\gwp^{(r)}(\mu_{c,\varepsilon},\nu_{c,\varepsilon}) \leq \gwp^{(r)}(\mu,\nu) + 4C\varepsilon$.
\end{proof}

The proof also requires a reverse process which we call {\em extension}.  Recall that $\cX \times \cY$ admits the following partition:
\begin{equation*}
\cX \times \cY = \bigsqcup_{i \in [k],\; j\in [l]}(\cX_i \times \cY_j).
\end{equation*}  Given an atomic coupling measure $\overline{\pi} \in \Pi(\mu_{c,\varepsilon},\nu_{c,\varepsilon})$, we define its $\varepsilon$-{\em extension} by 
\begin{equation*}
\overline{\pi}_{e,\varepsilon}= \sum_{i\in [k],\;j\in[l]}\delta_{\cX_i \times \cY_j}\lambda_{ij}\mu \otimes \nu, \quad \lambda_{ij}= \begin{cases}
0 & \text{if } \mu(\cX_i)\nu(\cY_j)=0\\
\dfrac{\overline{\pi}((\bx_i,\by_j))}{\mu(\cX_i)\nu(\cY_j)} & \text{otherwise}
\end{cases}.
\end{equation*}

\begin{lemma}\label{thm:extensionbound}
We have the following inequality: 
\begin{equation}\label{eq: concentration inequality}
\gwp(\mu,\nu) \leq \gwp(\mu_{c,\varepsilon},\nu_{c,\varepsilon}) + 4C\varepsilon.
\end{equation}
\end{lemma}

\begin{proof}[Proof of Lemma \ref{thm:extensionbound}]
We claim that the extension satisfies $\overline{\pi}_{e,\varepsilon} \in \Pi(\mu,\nu)$.  In particular, it suffices to prove that, for any Borel $A \subset \cX_i$ where $i \in [k]$ and any Borel $B \subset \cY_j$ where $j \in [l]$, one has the following
\begin{equation*}
\overline{\pi}_{e,\varepsilon}(A \times \cY) = \mu(A), \quad \mathrm{and} \quad \overline{\pi}_{e,\varepsilon}(\cX \times B) = \nu(B).
\end{equation*}
Indeed, by the definition of $\overline{\pi}_{e,\varepsilon}$, we have
\begin{equation*}
\overline{\pi}_{e,\varepsilon}(A \times \cY)= \sum_{j \in [l]}\delta_{\cX_i \times \cY_j}\lambda_{ij}\mu\times \nu(A \times \cY) 
= \sum_{j \in [l]}\lambda_{ij}\mu(A)\nu(\cY_j) =\mu(A)\sum_{j \in [l]}\lambda_{ij}\nu(\cY_j).
\end{equation*}
Suppose $\mu(\cX_i) =0$.  Then $\lambda_{ij}=0$ and $\mu(A)=0$.  Suppose on the other hand $\mu(\cX_i) \neq 0$.  Then
\begin{equation*}
\overline{\pi}_{e,\varepsilon}(A \times \cY)= \frac{\mu(A)}{\mu(\cX_i)}\sum_{j \in [l]:\; \nu(\cY_j) \neq 0}\overline{\pi}(\bx_i,\by_j)=\frac{\mu(A)}{\mu(\cX_i)} \mu_{c,\varepsilon}(\bx_i) =\mu(A).
\end{equation*}
In either case, we see that $\mu(A)=\overline{\pi}_{e,\varepsilon}(A \times \cY)$.  Similarly, we also have $\overline{\pi}_{e,\varepsilon}(\cX \times B) = \nu(B)$.  Therefore we conclude $\overline{\pi}_{e,\varepsilon} \in \Pi(\mu,\nu)$.  By \eqref{eq: lipszhitz inequality} we have 
\begin{equation*}\label{eq:concentration_ineq2}
\int_{(\cX \times \cY)^{2}}c(\overline{\bx}_1,\overline{\by}_1,\overline{\bx}_2,\overline{\by}_2)~d\overline{\pi} \leq \int_{(\cX \times \cY)^{2}}c(\overline{\bx}_1,\overline{\by}_1,\overline{\bx}_2,\overline{\by}_2)~d\pi_{e,\varepsilon} + 4C\varepsilon.
\end{equation*}
Then by minimizing over $\overline{\pi} \in \Pi(\mu_{c,\varepsilon},\nu_{c,\varepsilon})$ we conclude \eqref{eq: concentration inequality}.
\end{proof}

\begin{lemma}\label{thm:extensionhierarchy}
For any $\overline{\bP} \in \Pi_r(\mu_{c,\varepsilon},\nu_{c,\varepsilon})$, we define an $\varepsilon-$extension of $\overline{\bP}$ by 
\begin{equation*}
\overline{\bP}_{e,\varepsilon}= \sum_{\substack{i_1,\dots,i_{2r}\in [k]\\j_1,\dots,j_{2r}\in[l]}}\delta_{\prod_{s=1}^{2r}(\cX_{i_s} \times \cY_{j_s})}\lambda_{\substack{i_1,\dots,i_{2r}\\j_1,\dots,j_{2r}}}\mu^{2r} \otimes \nu^{2r},
\end{equation*}
where
\begin{equation*}
\lambda_{ij}= \begin{cases}
0 & \text{if } \prod_{s=1}^{2r}\mu(\cX_{i_s})\nu(\cY_{j_s})=0\\
\dfrac{\overline{\bP}((\bx_{i_1},\by_{j_1},\dots,\bx_{i_{2r}},\by_{j_{2r}}))}{\prod_{s=1}^{2r}\mu(\cX_{i_s})\nu(\cY_{j_s})} & \text{otherwise}
\end{cases}.
\end{equation*}
Then $\overline{\bP}_{e,\varepsilon} \in \Pi_r(\mu,\nu)$. Consequently, we have the following bound.
\begin{equation}\label{ineq 4}
\gwp^{(r)}(\mu,\nu) \leq \gwp^{(r)}(\mu_{c,\varepsilon},\nu_{c,\varepsilon}) + 4C\varepsilon.
\end{equation}
\end{lemma}

\begin{proof}
First, $\overline{\bP}_{e,\varepsilon}$ inherits \eqref{eq:cond_sym} from $\overline{\bP}$.  Second, we can verify \eqref{eq:cond_mar} in a similar fashion to the proof of Lemma \ref{thm:concentrationhierarchy}.  It remains to prove $\overline{\bP}_{e,\varepsilon}$ is positive semidefinite.  

Let $0 \leq t \leq r$, $f \in \cB((\cX \times \cY)^{2t})$, and $g \in \cB_+((\cX \times \cY)^{2r-2t})$.  Define the following step functions 
\begin{equation*}
\begin{aligned}
\overline{f}(\overline{\bx}_1,\overline{\by}_1,\dots,\overline{\bx}_t,\overline{\by}_t) :=& \sum_{\substack{i_1,\dots,i_t \in [k]\\ j_1,\dots,j_t \in [l]}}\int_{\prod_{s=1}^{t}(\cX_{i_s} \times \cY_{j_s})}f(\overline{\bx}_1,\overline{\by}_1,\dots,\overline{\bx}_t,\overline{\by}_t)~d \mu^{t} \otimes \nu^{t} \delta_{\prod_{s=1}^t(\cX_{i_s}\times \cY_{j_s})}\\
\overline{g}(\overline{\bx}_1,\overline{\by}_1,\dots,\overline{\bx}_{2r-2t},\overline{\by}_{2r-2t}) := & \sum_{\substack{i_1,\dots,i_{2r-2t} \in [k]\\ j_1,\dots,j_{2r-2t} \in [l]}}\int_{\prod_{s=1}^{2r-2t}(\cX_{i_s} \times \cY_{j_s})}g(\bx_{i_1},\by_{j_1},\dots,\bx_{i_{2r-2t}},\by_{j_{2r-2t}}) \\
& \qquad\qquad\qquad\qquad d\mu^{2r-2t}\otimes \nu^{2r-2t}\delta_{\prod_{s=1}^{2r-2t}(\cX_{i_s}\times \cY_{j_s})}.
\end{aligned}
\end{equation*}
Then, by the definition of $\bP_{e,\varepsilon}$, we have
\begin{align*}
&\int_{(\cX \times \cY)^{2r}}f(\overline{\bx}_1,\overline{\by}_1,\dots,\overline{\bx}_t,\overline{\by}_t)f(\overline{\bx}_{t+1},\overline{\by}_{t+1},\dots,\overline{\bx}_{2t},\overline{\by}_{2t})g(\overline{\bx}_{2t+1},\overline{\by}_{2t+1},\dots,\overline{\bx}_{2r},\overline{\by}_{2r})~d\overline{\bP}_{e,\varepsilon}\\
=\;& \int_{(\cX \times \cY)^{2r}}\overline{f}(\overline{\bx}_1,\overline{\by}_1,\dots,\overline{\bx}_t,\overline{\by}_t)\overline{f}(\overline{\bx}_{t+1},\overline{\by}_{t+1},\dots,\overline{\bx}_{2t},\overline{\by}_{2t})\overline{g}(\overline{\bx}_{2t+1},\overline{\by}_{2t+1},\dots,\overline{\bx}_{2r},\overline{\by}_{2r})~d\overline{\bP} \geq 0.
\end{align*}
The last inequality comes from the positive semidefiniteness of $\overline{\bP}$.  As such, we have $\overline{\bP}_{e,\varepsilon} \in \Pi_r(\mu,\nu)$.  By applying \eqref{eq: lipszhitz inequality}, we obtain
\begin{equation*}
\int_{(\cX \times \cY)^{2r}}c(\overline{\bx}_1,\overline{\by}_1,\overline{\bx}_2,\overline{\by}_2)~d\overline{\bP}\leq \int_{(\cX \times \cY)^{2r}}c(\overline{\bx}_1,\overline{\by}_1,\overline{\bx}_2,\overline{\by}_2)~d\overline{\bP}_{c,\varepsilon} + 4C\varepsilon.
\end{equation*}
Then, by minimizing over $\overline{\bP} \in \Pi_r(\mu_{c,\varepsilon},\nu_{c,\varepsilon})$, we conclude \eqref{ineq 4}.
\end{proof}

\begin{proof}[Proof of Theorem \ref{thm: convergence}]
By combining Lemma \ref{thm:concentrationbound} and \ref{thm:extensionbound} we obtain 
\begin{equation*}
|\gwp (\mu_{c,\varepsilon},\nu_{c,\varepsilon}) - \gwp (\mu,\nu) | \leq 4 C \varepsilon.  
\end{equation*}
Similarly, by combining Lemma \ref{thm:concentrationhierarchy} and \ref{thm:extensionhierarchy}, we obtain 
\begin{equation*}
|\gwp^{(r)}(\mu_{c,\varepsilon},\nu_{c,\varepsilon}) - \gwp^{(r)}(\mu,\nu) | \leq 4 C \varepsilon.
\end{equation*}
This implies
\begin{equation*}
\begin{aligned}
&\left|\gwp(\mu,\nu)- \gwp^{(r)}(\mu,\nu) \right|\\
\leq\;& \left|\gwp(\mu,\nu)- \gwp(\mu_{c,\varepsilon},\nu_{c,\varepsilon}) \right| + \left|\gwp(\mu_{c,\varepsilon},\nu_{c,\varepsilon})- \gwp^{(r)}(\mu_{c,\varepsilon},\nu_{c,\varepsilon}) \right|+ \left|\gwp^{(r)}(\mu_{c,\varepsilon},\nu_{c,\varepsilon})-\gwp^{(r)}(\mu,\nu) \right|\\
\leq \;& 8C \varepsilon + \left|\gwp(\mu_{c,\varepsilon},\nu_{c,\varepsilon})- \gwp^{(r)}(\mu_{c,\varepsilon},\nu_{c,\varepsilon}) \right|.
\end{aligned}
\end{equation*}
Here, $\mu_{c,\varepsilon}$ and $\nu_{c,\varepsilon}$ are discrete measures.  By \cite{tran:25}, we have $\gwp^{(r)}(\mu_{c,\varepsilon},\nu_{c,\varepsilon}) \rightarrow \gwp(\mu_{c,\varepsilon},\nu_{c,\varepsilon})$ as $r \rightarrow \infty$.  As such, we have
\begin{equation*}
\limsup_{r \to \infty}\left|\gwp(\mu,\nu)- \gwp^{(r)}(\mu,\nu) \right| \leq 8C\varepsilon \quad \text{for all} \quad \varepsilon >0,
\end{equation*}
and hence  $\lim_{r\to \infty}\gwp^{(r)}(\mu,\nu)= \gwp(\mu,\nu)$.
\end{proof}

\section{Metric Properties} \label{sec:metric}


Let $\FX_{p,q}$ denote the set of compact metric measure spaces. 
Given a pair of metric measure spaces $\mathbb{X}:=(\CX,d_{\CX},\mu)$ and $\mathbb{Y}:=(\CY,d_{\CY},\nu)$, we define the distance measure
\begin{equation*}
\gw_{p,q}^{(r)}(\mathbb{X},\mathbb{Y})= \left[ \gwp^{(r)}(\mathbb{X},\mathbb{Y})\right]^{1/p},
\end{equation*}
where the cost function $c$ is chosen as in \eqref{eq:lp_definition}.  The main result of this section is to show that $\gw_{p,q}^{(r)}$ specifies a pseudo-metric in a similar fashion the $L_{p,q}$-distortion distance $\gw_{p,q}$ does over the space of metric measure spaces.  This is the formal result: 
\begin{theorem}\label{thm: metric}
Let $r$ be a positive integer.  Then $\gw := \gw_{p,q}^{(r)}$ is a pseudo-distance over $\FX_{p,q}$; that is, it satisfies
\begin{enumerate}
\item $\gw (\bbX,\bbX) =0$,
\item (Non-negativity) $\gw (\bbX,\bbY) \geq 0 $ for all $\bbX, \bbY \in \FX_{p,q}$,
\item (Symmetry) $\gw (\bbX,\bbY)= \gw (\bbY,\bbX)$ for all $\bbX, \bbY \in \FX_{p,q}$, and
\item (Triangle Inequality) $\gw (\bbX,\bbZ) \leq \gw (\bbX,\bbY) + \gw (\bbY,\bbZ)$ for all $\bbX, \bbY \in \FX_{p,q}$.
\end{enumerate}
\end{theorem}
Among these properties, non-negativity and symmetry are relatively straightforward to prove, and we omit the proof of these steps.  The main focus of this section is to prove the triangle inequality.  Key to this is an intermediate result frequently referred to as the {\em gluing lemma} that glues together transportation plans with the appropriate marginals.  We begin by stating the precise result, before providing some historical context.

\begin{lemma}[Gluing Lemma] \label{thm:gluing}
Let $\bbX =(\cX,d_{\cX},\mu)$, $\bbY =(\cY,d_{\cY},\nu)$, and $\bbZ =(\cZ,d_{\cZ},\kappa)$ be metric measure spaces in $\FX_{p,q}$.  Let $\bP \in \Pi_r(\mu,\nu)$ and $\bQ \in \Pi_r(\nu,\kappa)$.  Then there exists $\bS \in \cP((\cX \times \cY \times \cZ)^{2r})$ such that (i) $\bS|_{(\cX \times \cY)^{2r}}= \bP$ and $\bS|_{(\cY \times \cZ)^{2r}}= \bQ$; and, moreover, (ii) if we define $\bR:= \bS|_{(\cX \times \cZ)^{2r}}$, then we also have $\bR \in \Pi_r(\mu,\kappa)$.
\end{lemma}

Lemma \ref{thm:gluing} generalizes, or can be viewed as the analog of, the gluing lemma in a number of simpler contexts.  The simplest instance of the gluing lemma appears in the context of classical optimal transport; see, for instance, \cite{Vil:08}.  The work in \cite{Sturm:23} proves an analogous result for the GW problem.  As such, the main question we ask in this paper is whether a similar type of result also holds for the proposed moment-SOS hierarchy for the GW problem.  First, we note that the work in \cite{tran:25} proves such a result in settings where the spaces $\CX$, $\CY$, and $\CZ$ are discrete -- see, in particular, \cite[Lemma 5.2]{tran:25}.  As such, Lemma \ref{thm:gluing} can be viewed as the extension of \cite[Lemma 5.2]{tran:25} for continuous distributions.  

We make some important remarks:  While the statement of Lemma \ref{thm:gluing} may superficially appear to be the continuous extension of \cite[Lemma 5.2]{tran:25}, there are important underlying conceptual differences.  Specifically, the analysis in Lemma \ref{thm:gluing} deals directly with probability distributions.  On the other hand, the analysis in \cite[Lemma 5.2]{tran:25} deals directly with the optimization formulation \eqref{eq:schmudgen_dgw} without appealing to any probabilistic interpretation.  In a sense, our analysis uncovers a hidden probabilistic interpretation that underlies the analyses in \cite[Lemma 5.2]{tran:25}.

\begin{proof}[Proof of Lemma \ref{thm:gluing}]
The key step of the proof is the construction of $\bS$.  Given measurable subsets $A \subset \cX^{2r}$, $B \subset \cY^{2r}$, and $C \subset \cZ^{2r}$, we define  
\begin{equation*}
\bS(A \times B \times C) = 
\begin{cases}
0 & \text{if }\nu^{2r}(B)=0\\
\frac{\bP(A \times B)\bQ(B \times C)}{\nu^{2r}(B)} & \text{otherwise}
\end{cases}.
\end{equation*}
Note that $\bS$ is well-defined because $\cB(\cX^{2r}) \times \cB(\cY^{2r}) \times \cB(\cZ^{2r})$ generates $\cB((\cX \times \cY \times \cZ)^{2r})$.  It is straightforward to check the relations $\bS|_{(\cX \times \cY)^{2r}}= \bP$ and $\bS|_{(\cY \times \cZ)^{2r}}= \bQ$, and we omit these steps.  The last check to perform is that $\bR:= S|_{(\cX \times \cZ)^{2r}}$ belongs to $\Pi_r(\mu,\kappa)$.  To do this, there are three properties we need to check.

[\eqref{eq:cond_sym}]: Let $\sigma \in S_{2r}$ be a permutation, and let $X \subset \CX^{2r}$ and $Z \subset \CZ^{2r}$ be measurable subsets.  We then have
\begin{equation*}
\begin{aligned}
&\bR(\overline{\sigma}(X \times Z))= \bS(\overline{\sigma}(X \times \CY^{2r} \times Z))= \frac{\bP(\overline{\sigma}(X \times \CY^{2r}))\bQ(\overline{\sigma}(\CY^{2r} \times Z))}{\nu^{2r}(\overline{\sigma}(\CY^{2r})})\\
\overset{(a)}{=}\;& \frac{\bP(X \times \cY^{2r})\bQ(B \times \CY^{2r})}{\nu^{2r}(\CY^{2r})}=\bS (X \times \CY^{2r} \times Z)= \bR (X \times Z).
\end{aligned}
\end{equation*}
Here, (a) uses the assumption $\bP$ and $\bQ$ satisfy \eqref{eq:cond_sym}.

[\eqref{eq:cond_mar}]:  Let $X \in \CX$, $Z \in \CZ$, $X' \subset \CX^{2r-1}$, and $Z' \subset \CZ^{2r-1}$ be measurable subsets.  We then have
        \begin{align*}
            &\bR(X \times X' \times Z') = \bS(X \times \cY^{2r} \times X' \times Z') = \frac{\bP(X \times X' \times \cY^{2r})\bQ(\cY^{2r} \times Z')}{\nu^{2r}(\cY^{2r})}\\
            \overset{(b)}{=}\;& \mu(X)\frac{\bP(X' \times \cY^{2r})\bQ(\cY^{2r}\times Z')}{\nu^{2r}(\cY^{2r})}= \mu(X)\bS(X' \times \cY^{2r} \times Z')= \mu(X)\bR(X' \times Z').
        \end{align*}
        Here, (b) uses the assumption $\bP$ satisfies \eqref{eq:cond_mar}.  Using a similar sequence of steps, we also have $\bR(Y \times X' \times Z') = \nu(Y) \times \bR(X' \times Z')$.  This verifies the marginal conditions.

[\eqref{cond: PSD for cts}]:  This part of the proof appears long but is conceptually quite straightforward.  In summary, verifying \eqref{cond: PSD for cts} requires us to establish an inequality.  To do so, we will rely on the concentration argument we introduced in Section \ref{sec: moment for cts GW} to pass over to a suitable discretized problem, and on which we apply the gluing lemma in \cite[Lemma 5.2]{tran:25}.

{\em Step 1 (Discretization):}  Recall that $(\CX,d_{\CX})$, $(\CY,d_{\CY})$, and $(\CZ,d_{\CZ})$ are compact metric spaces, and hence totally bounded.  As such, there exists partitions
\begin{equation*}
\cX = \bigsqcup_{i=1}^k\cX_i, \quad  \cY = \bigsqcup_{j=1}^l\cY_i \quad \mathrm{and} \quad \cZ = \bigsqcup_{v=1}^u\cZ_v
\end{equation*}
such that for any $i \in [k], j \in [l], v \in [u]$, there exists grid points $\bx_i \in \cX _i$, $\by_j \in \cY_j$ and $\bz_v \in \cZ_v$ so that $\cX_i \subset B(\bx_i,1/s) \; \forall i \in [k]$, $\cY_j \subset B(\by_j,1/s)\; \forall j \in [l]$, and $\cZ_v \subset B(\bz_v,1/s)\; \forall v \in [u]$.  We denote 
\begin{equation*}
\mathbf{Point}[t,s](\bbX,\bbY) := \left\{(\bx_{i_1},\by_{j_1},\dots,\bx_{i_t},\by_{j_t})\;:\;i_1,\dots,i_t \in [k],\; j_1,\dots,j_t \in [l]\right\}.
\end{equation*}
We define the following collection of step functions 
\begin{equation*}
\mathbf{SF}[t,s](\bbX,\bbY)\:= \left\{\sum_{\substack{i_1,\dots,i_t \in [k]\\ j_1,\dots,j_t \in [l]}}\lambda_{\substack{i_1,\dots,i_t \in [k]\\ j_1,\dots,j_t \in [l]}}\delta_{\prod_{w=1}^s(\cX_{i_w}\times \cY_{j_w})}\;:\;\lambda_{\substack{i_1,\dots,i_s \in [k]\\ j_1,\dots,j_s \in [l]}} \in \RR \right\}.
\end{equation*}
Last, we denote
\begin{equation*}
\mathbf{PSF}[t,s](\bbX,\bbY):= \left\{f \in \mathbf{SF}(s,t)\;:\; f \geq 0\right\}.
\end{equation*}
We then note the following:
\begin{itemize}
\item $\bigcup_{t=1}^{\infty}\mathbf{Point}[t,s](\bbX,\bbY)$ is dense in $(\cX \times \cY)^{t}$, and
\item $\bigcup_{t=1}^{\infty}\mathbf{SF}[t,s](\bbX,\bbY)$ is dense in $\cB((\cX \times \cY)^t)$ with respect to the $L^{\infty}$ topology,
\item $\bigcup_{t=1}^{\infty}\mathbf{PSF}[t,s](\bbX,\bbY)$ is dense in $\cB_+((\cX \times \cY)^t)$ with respect to the $L^{\infty}$ topology.
\end{itemize}
In particular, to check that \eqref{cond: PSD for cts} holds, it suffices to check that
\begin{multline}\label{cond: discrete PSD}
        \int_{(\cX \times \cY)^{2r}}f(\bx_1,\by_1,\dots,\bx_t,\by_{t})f(\bx_{t+1},\by_{t+1},\dots,\bx_{2t},\by_{2t})g(\bx_{2t+1},\by_{2t+1},\dots,\bx_{2r},\by_{2r})~d\bP \geq 0 \\ \forall \;s \in \NN_{>0},\; 0 \leq t \leq r,\; f \in \mathbf{SF}(t,s),\; g \in \mathbf{PSF}(t,s).
\end{multline}

{\em Step 2:}  We proceed to \eqref{cond: discrete PSD} for fixed $s$.
    To do that, we recall the $1/s-$concentration of $\bP,\bQ$ and $\bR$ in Section \ref{sec: moment for cts GW}.  By Lemma \ref{thm:concentrationhierarchy}, we have  
    \begin{displaymath}
        \bP_{c,1/s} \in \Pi_r(\mu_{c,1/s},\nu_{c,1/s}), \quad \text{and} \quad \bQ_{c,1/s} \in \Pi_r(\nu_{c,1/s},\kappa_{c,1/s}).
    \end{displaymath}
    and 
    \begin{align*}
        &\bR_{c,1/s}((\bz_{v_1},\bx_{i_1},\dots,\bz_{v_{2r}},\bx_{i_{2r}}))\\
        =\;& \begin{cases}
            0 & \text{if } \prod_{w=1}^{2r}\nu_{c,1/s}(\by_{j_w}) =0\\
            \dfrac{\bP_{c,1/s}((\bx_{i_1},\by_{j_1},\dots,\bx_{i_{2r}},\by_{j_{2r}}))\bQ_{c,1/s}((\by_{j_1},\bz_{v_1},\dots,\by_{j_{2r}},\bz_{v_{2r}}))}{\prod_{w=1}^{2r}\nu_{c,1/s}(\by_{j_w})}& \text{otherwise}
        \end{cases}.
    \end{align*}

Given $f \in \mathbf{SF}(t,s)$ and $g \in \mathbf{PSF}(t,s)$, we have
\begin{equation*}
\begin{aligned}
&\int_{(\cX \times \cZ)^{2r}}f(\overline{\bx}_1,\overline{\bz}_1,\dots,\overline{\bx}_t,\overline{\bz}_{t})f(\overline{\bx}_{t+1},\overline{\bz}_{t+1},\dots,\overline{\bx}_{2t},\overline{\bz}_{2t})g(\overline{\bx}_{2t+1},\overline{\bz}_{2t+1},\dots,\overline{\bx}_{2r},\overline{\bz}_{2r})~d\bR\\
=\;&\int_{(\cX \times \cZ)^{2r}}f(\overline{\bx}_1,\overline{\bz}_1,\dots,\overline{\bx}_t,\overline{\bz}_{t})f(\overline{\bx}_{t+1},\overline{\bz}_{t+1},\dots,\overline{\bx}_{2t},\overline{\bz}_{2t})g(\overline{\bx}_{2t+1},\overline{\bz}_{2t+1},\dots,\overline{\bx}_{2r},\overline{\bz}_{2r})~d\bR_{c,1/s}.
\end{aligned}
\end{equation*}
Then, from Lemma 5.2 (see e.g.,\cite{tran:25}), $\bR_{c,1/s} \in \Pi_r(\kappa_{c,1/s},\mu_{c,1/s})$.  In particular, this means that $\bR_{c,1/s}$ satisfies \eqref{cond: PSD for cts}, and hence the above expression is non-negative.  Now take limits $s \rightarrow 0$ to conclude that $\bR$ satisfies \eqref{cond: PSD for cts}.  This completes the proof.
\end{proof}


We now have all the ingredients to prove the triangular inequality.

\begin{lemma}\label{thm:triangle}
Let $\bbX$, $\bbY$ and $\bbZ$ be metric measure spaces.  We then have the following
\begin{equation*}
\gw (\bbX,\bbZ) \leq \gw (\bbX,\bbY) + \gw (\bbY,\bbZ).
\end{equation*}
\end{lemma}

\begin{proof}[Proof of Theorem \ref{thm:triangle}]
For any $\varepsilon >0$, we choose $\bP \in \Pi_r(\mu,\nu)$ and $\bQ \in \Pi_r(\nu,\kappa)$ such that 
\begin{equation*}
\begin{aligned}
\Big|\gw (\bbX,\bbY) - \Big(\int_{(\cX \times \cY)^{2r}}| d_{\cX} (\bx_1,\bx_2)^q - d_{\cY} (\by_1,\by_2)^q |^p~d\bP \Big)^{1/p}\Big| \leq & \varepsilon,\\
\Big|\gw (\bbY,\bbZ) - \Big(\int_{(\cY \times \cZ)^{2r}}| d_{\cY} (\by_1,\by_2)^q - d_{\cZ} (\bz_1,\bz_2)^q |^p~d\bQ \Big)^{1/p}\Big| \leq & \varepsilon.
\end{aligned}
\end{equation*}
    We construct $\bS \in \cP((\cX \times \cY \times \cZ)^{2r})$ and $\bR \in \Pi_r(\kappa,\mu)$ as in the proof of the Gluing Lemma \ref{thm:gluing}. By the triangular inequality in $L^p$, we evaluate
\begin{equation*}
\begin{aligned}
& \big(\int_{(\cX \times \cZ)^{2r}}| d_{\cX} (\bx_1,\bx_2)^q - d_{\cZ} (\bz_1,\bz_2)^q |^p~d\bR \big)^{1/p} \\
=~ & \big(\int_{(\cX \times \cY \times \cZ)^{2r}}| d_{\cX} (\bx_1,\bx_2)^q - d_{\cZ} (\bz_1,\bz_2)^q |^p~d\bS \big)^{1/p} \\
\leq~ & \big(\int_{(\cX \times \cY \times \cZ)^{2r}}| d_{\cX} (\bx_1,\bx_2)^q - d_{\cY} (\by_1,\by_2)^q |^p~d\bS \big)^{1/p} \\
& \qquad \qquad +~ \big(\int_{(\cX \times \cY \times \cZ)^{2r}}| d_{\cY} (\by_1,\by_2)^q - d_{\cZ} (\bz_1,\bz_2)^q |^p~d\bS \big)^{1/p}.
\end{aligned}
\end{equation*}

This implies
\begin{equation*}
\gw (\bbX,\bbZ) \leq \gw (\bbX,\bbY) + \gw (\bbY,\bbZ) +2\varepsilon .
\end{equation*}
Since $\epsilon$ is arbitrary, we conclude $\gw (\bbX,\bbZ) \leq \gw (\bbX,\bbY) + \gw (\bbY,\bbZ)$.
\end{proof}

Last, we note that Lemma \ref{thm:triangle} also completes the proof for the Theorem \ref{thm: metric}.

\section{Sample Complexity} \label{sec:consistency}

The goal of the prior sections was to introduce the moment-SOS hierarchy $\gwp^{(r)}$ as an approximation of the desired $\gwp$ distance, as well as to understand the analytical properties underlying $\gwp^{(r)}$.  In this section, the basic question is one concerning statistical consistency.  Consider metric measure spaces $\bbX=(\cX,d_{\cX},\mu)$ and $\bbY=(\cY,d_{\cY},\nu)$.  Suppose one draws samples from $\mu$ and $\nu$, using which one specifies the empirical distributions $\hat{\mu}_n$ and $\hat{\nu}_n$ respectively.  Suppose one proceeds to compute the discrete GW distance between $\hat{\mu}_n$ and $\hat{\nu}_n$ via the moment-SOS hierarchy \eqref{eq:schmudgen_dgw}.  How does the resulting distance compare with the GW distance (as well as the SOS analogs) between $\mu$ and $\nu$?

In Section \ref{sec: moment for cts GW} we gave a partially-related answer to this question in which we discretize the underlying metric space $\CX$ and $\CY$ uniformly using $\epsilon$-balls.  The treatment in this section is very different because the samples we draw from $\mu$ and $\nu$ are independent and identically distributed (i.i.d.), and as such, we do not get to exert control over how uniformly spread out these samples are.  The set-up we consider in this section is more typical in empirical statistics.

In what follows, we let $\hat{\mu}_n = \sum_{i=1}^n \delta_{X_i}, X_i \sim \mu$ denote the empirical measures of $\mu$, and likewise for $\hat{\nu}_n$.  The samples are i.i.d. draws from $\mu$ and $\nu$ respectively.  Our goal is to show the following holds, for all integers $r$:
\begin{equation}\label{eq:sample-convergence}
    \bbE\left[\Delta^{(r)}((\cX,d_{\cX},\hat{\mu}_n),(\cY,d_{\cY},\hat{\nu}_n))\right] \rightarrow \Delta^{(r)} ((\cX,d_{\cX},\mu),(\cY,d_{\cY},\nu)) \quad \text{ as } \quad  n \rightarrow \infty.
\end{equation}

{\bf Related work.}  We discuss the results in this section in relation to relevant prior work.

The first piece of related work is that in \cite{ZGMS:24}.  Specifically, the authors in \cite{ZGMS:24} prove a result that is in essence \eqref{eq:sample-convergence}, but only in the case $p=q=2$ (see Theorem 2 of \cite{ZGMS:24} for the precise statement).  In contrast, our result applies to more general $L_{p,q,}$-distortion distances.  The reason the work in \cite{ZGMS:24} is restricted to the setting $p=q=2$ is that the analysis depends on being able to view the cost function $c$ as a polynomial, on which certain algebraic calculations are possible.  

The second piece of related work is that in \cite{WB:19}, which concerns the convergence of empirical measures to the underlying distribution in the Wasserstein distance $\mathbb{E}[ W(\hat{\mu}_n, \mu)] \rightarrow 0$.  There are some distinctions between the results in \cite{WB:19} and the results in this section.  First, the GW problem is a different problem from classical OT.  Second, the quantity $\mathbb{E}[ W(\hat{\mu}_n, \mu)]$ measures the difference between two {\em distributions}, whereas the statement in \eqref{eq:sample-convergence} concerns the {\em value} of the GW problem -- these are fundamentally different questions.  Nevertheless, the ideas in \cite{WB:19} are useful for us; in particular, in what follows, our analyses of using dyadic partitioning argument will closely follow that of \cite{WB:19}.

{\bf Simplification to a common space.}  To prove \eqref{eq:sample-convergence}, we will instead show the following:
\begin{equation} \label{eq:sampling_result}
\mathbb{E}\left[\gwp ((\cX,d_{\cX},\hat{\mu}_n),(\cX,d_{\cX},\mu))\right] \rightarrow 0 \quad \text{ as } \quad n \rightarrow \infty.
\end{equation}
The difference from \eqref{eq:sample-convergence} is that the underlying metric space in the source and target are identical in \eqref{eq:sampling_result}.

First, let us explain how \eqref{eq:sampling_result} leads to \eqref{eq:sample-convergence}.  By combining the triangle inequality in Theorem \ref{thm: metric} with Proposition \ref{thm:moment_lowerbounds} we have the following
\begin{equation*}    
\begin{aligned}
\Delta^{(r)} (\hat{\mu}_n,\hat{\nu}_n) & \leq \Delta^{(r)} (\hat{\mu}_n,\mu) + \Delta^{(r)} (\mu,\nu) + \Delta^{(r)} (\nu,\hat{\nu}_n) \\
& \leq \Delta (\hat{\mu}_n,\mu) + \Delta^{(r)} (\mu,\nu) + \Delta (\nu,\hat{\nu}_n).
\end{aligned}
\end{equation*}
In a similar fashion, we also have
\begin{equation*}    
\Delta^{(r)} (\mu,\nu) \leq \Delta (\mu,\hat{\mu}_n) + \Delta^{(r)} (\hat{\mu}_n,\hat{\nu}_n) + \Delta (\hat{\nu}_n,\nu).
\end{equation*}
Subsequently, by applying expectations on both sides (with respect to sampling the data-points), and by \eqref{eq:sampling_result} (as well as the corresponding statement for $\cY$) we conclude \eqref{eq:sample-convergence}.

The reason we make the simplification to \eqref{eq:sampling_result} is because, in settings where the metric spaces are equal, there are explicit constructions of transportion plans $\pi$ between $\hat{\mu}_n$ and $\hat{\nu}_n$ that evaluate to a small objective value in \eqref{eq:sampling_result}.

\subsection{Main result}

\begin{definition}[$\varepsilon$-dimension of a set]
The $\varepsilon$-dimension of the set $\cS \subset \cX$ is the quantity
\begin{equation*}
d_\varepsilon(\cS) := \dfrac{\log \mathcal N_\varepsilon (\cS)}{-\log \varepsilon}.
\end{equation*}
\end{definition}

\begin{definition}[Covering number of a set]
Given a measure $\mu$ on $\cX$, the $(\varepsilon, \tau)$-covering number of $\cX$ is the quantity
\begin{equation*}
\mathcal N (\mu, \tau) := \inf \{ \mathcal N_\varepsilon (\cS): \mu(\cS) \geq 1 - \tau\},
\end{equation*}
and the $(\varepsilon, \tau)$-dimension is
\begin{equation*}
d_\varepsilon (\mu, \tau) := \dfrac{\log \mathcal N_\varepsilon(\mu, \tau)}{-\log \varepsilon}.
\end{equation*}
We will also denote $\mathcal N_\varepsilon(\mu, 0) := \mathcal N_\varepsilon(\mu)$ and $d_\varepsilon (\mu, 0)= d_\varepsilon (\mu)$.
\end{definition}

The following is the main result of this section.

\begin{theorem} \label{thm:convergencerate}
Let $p \in [1, \infty)$.  Suppose there exists an $\varepsilon' \leq 1$ and $s > 2p$ such that
\begin{equation*}
d_\varepsilon(\mu, \varepsilon^{\frac{2s}{s-2p}}) \leq s
\end{equation*}
for all $\varepsilon \leq \varepsilon'$. Then
\begin{equation*}
\mathbb{E}[\gwp(\mu, \hat{\mu}_n)] \leq C_1n^{-pq/s} + \frac{3}{2} n^{-p/s} + C_2n^{-1/2},
\end{equation*}
where 
\begin{equation*}
C_1 := 3^{3pq} + \dfrac{3^{3pq}}{3^{s/2 - pq} - 1} + 3^{3\alpha + 1} \text{ and } C_2 :=  \frac{3}{2} \left(\frac{3}{\varepsilon'} \right)^{s/2}.
\end{equation*}
In particular, this means
\begin{equation*}
\mathbb{E}[\Delta(\mu, \hat{\mu}_n)] = \mathbb{E}[(\gwp(\mu, \hat{\mu}_n))^{1/p}] \lesssim n^{-1/s}.
\end{equation*}
\end{theorem}

To be clear, Theorem \ref{thm:convergencerate} is adapts Proposition 5 in \cite{WB:19} to our set-up (it is also half of the statement in Theorem 1 in \cite{WB:19}).  Theorem \ref{thm:convergencerate} also implies a convergence rate of $n^{-1/\max \{ d^{\star} (\mu), d^{\star} (\nu) \} }$, where we define
$$
d^{\star} (\mu) := \inf \{ s \in (2p,\infty) : \underset{\epsilon \rightarrow 0}{\lim \sup} ~
d_\varepsilon(\mu, \varepsilon^{\frac{2s}{s-2p}}) \leq s
\}.
$$
The quantity $d^{\star}$ captures some notion of an intrinsic dimension of the probability distribution -- a more detailed discussion can be found in \cite{WB:19}.

\subsection{Proof of Theorem \ref{thm:convergencerate}}

The proof of Theorem \ref{thm:convergencerate} follows closely to the results in \cite{WB:19}.

\begin{definition}[\cite{WB:19}]
A {\em dyadic partition} of a set $\cS \subseteq \cX$ with parameter $\delta<1$ is a sequence of partitions $\{ \mathfrak{Q}^k \}$ satisfying
\begin{enumerate}
\item Each $\mathfrak{Q}^k$ is a partition of $\cS$.
\item Every $\cQ \in \mathfrak{Q}^k$ is Borel, and satisfies $\mathrm{diam}(\cQ) \leq \delta^k$.
\item If $\cQ^k \in \mathfrak{Q}^{k}$ and $\cQ^{k+1} \in \mathfrak{Q}^{k+1}$, then it must be that $\cQ^{k+1} \subset \cQ^k$ or $\cQ^{k+1} \cap \cQ^{k} = \emptyset$.
\end{enumerate}
\end{definition}

We introduce the following definition to simplify some of the calculations later on.  Given $\varepsilon>0$, we define the function $h(\varepsilon)$ as the largest cost that can be attained by pairs of points, each coming from a set with small diameter:
$$
h(\varepsilon) := \sup_{\cQ_1 , \cQ_2 : \mathrm{diam}(\cQ_i) \leq \varepsilon} \Big ( \sup_{\bx_1,\by_1 \in \cQ_1 , \bx_2,\by_2 \in \cQ_2} c (\bx_1,\by_1,\bx_2,\by_2) \Big ).
$$
For instance, if we consider taking $c$ to be the $L^{p,q}$ distortion distance between pairs of metric measure spaces, then $h(\varepsilon) \leq \varepsilon^{pq}$.

\begin{proposition} \label{thm:analog_proposition1}
Let $\lambda$ and $\mu$ be Borel probability measures on $\cX$.  Let $\cS \subset \cX$ be such that $\lambda(\cS)=\mu(\cS)$.  Suppose $\{ \mathfrak{Q}^{k} \}_{k=1}^{k^*}$ is a dyadic partition with parameter $\delta$.  Then
\begin{equation*}
\gwp(\lambda,\mu) \leq 3 \lambda(\cS) \times \sum_{k=0}^{k^*-1} h(\delta^{k-1}) \Big (\sum_{\cQ \in \mathfrak{Q}^{k}} | \lambda(\cQ) - \mu(\cQ)| \Big) + h(\delta^k) \lambda(\cS)^2.
\end{equation*}
\end{proposition}


Let $\lambda$ and $\mu$ be Borel probability measures on $X$.  Let $\cS \subset $ be such that $\lambda(\cS)=\mu(\cS)$, and let $\pi \in \Pi (\lambda,\mu)$ be a transportation plan on $\cS$.  Then suppose $\mathfrak{Q}$ is a partition of $\cS$ such that $\mathrm{diam}(\cQ) \leq \delta$ for all $\cQ \in \mathfrak{Q}$.  For such a partition, we define the terms
$$
t_{ijkl} := \int_{\bx_1 \in \cQ_i, \by_1 \in \cQ_j, \bx_2 \in \cQ_k, \by_2 \in \cQ_l} d \pi(\bx_1,\by_1) d \pi(\bx_2,\by_2).
$$
The proof of Proposition \ref{thm:analog_proposition1} proceeds by constructing an explicit transportation plan that moves mass in a way that respects the partitions $\mathfrak{Q}$, and iteratively through the refinements.  In doing so, we need to track how the terms of the form $t_{ijkl}$ contribute to the GW objective.

\begin{lemma}
One has 
$$
\sum_{i,j,k,l : i \neq j, k \neq l} t_{ijkl} \leq \Big ( \sum_i |\lambda(\cQ_i) - \mu(\cQ_i)| \Big)^2.
$$
\end{lemma}

\begin{proof}
First, notice that the terms in $i,j$ and the terms $k,l$ can be separated in the sum.  Hence one has
$$
\begin{aligned}
\sum_{i,j,k,l : i \neq j, k \neq l} t_{ijkl} = & \sum_{i,j,k,l : i \neq j, k \neq l} \int_{\bx_1 \in \cQ_i, \by_1 \in \cQ_j, \bx_2 \in \cQ_k, \by_2 \in \cQ_l} d \pi (\bx_1,\by_1) d \pi (\bx_2,\by_2) \\
= & \Big( \sum_{i,j : i \neq j} \int_{\bx_1 \in \cQ_i, \by_1 \in \cQ_j} d \pi (\bx_1,\by_1) \Big) \Big( \sum_{k,l : k \neq l} \int_{\bx_2 \in \cQ_k, \by_2 \in \cQ_l} d \pi (\bx_2,\by_2) \Big).\end{aligned}
$$
Consider $\sum_{i,j : i \neq j} \int_{x_1 \in \cQ_i, y_1 \in \cQ_j} d \pi (\bx_1,\by_1)$.  We fix a set $\cQ_i$ and consider the sum over the indices $j$.  One has
$$
\sum_{j : i \neq j} \int_{\bx_1 \in \cQ_i, \by_1 \in \cQ_j} d \pi (\bx_1,\by_1) = | \lambda (\cQ_i) - \mu(\cQ_i)|.
$$
This is because the integral accounts for all probability mass entering or exiting $Q_i$.  Now we sum over the indices $i$ to obtain the result.
\end{proof}

\begin{lemma}
One has 
$$
\sum_{i,j,k : i \neq j} t_{ijkk} \leq \lambda (\cS) \times \Big ( \sum_i |\lambda(\cQ_i) - \mu(\cQ_i)| \Big).
$$
\end{lemma}


\begin{proof}
The proof is similar in some parts.  In the first step, we decouple the sums in $i,j$ and the sums in $k$ to obtain
$$
\sum_{i,j,k : i \neq j} t_{ijkk} = \Big( \sum_{i,j : i \neq j} \int_{x_1 \in \cQ_i, y_1 \in \cQ_j} d \pi (\bx_1,y_1) \Big) \Big( \sum_{k} \int_{\bx_2 \in \cQ_k, \by_2 \in \cQ_k} d \pi (\bx_2,\by_2) \Big).
$$
In the previous result, we show that the term in the first parenthesis is bounded above by $\sum_i | \lambda (\cQ_i) - \mu(\cQ_i)|$.  The second sum $\sum_{k} \int_{\bx_2 \in \cQ_k, \by_2 \in \cQ_k} d \pi (\bx_2,\by_2)$ is (trivially) bounded above by $\lambda(\cS)$, which completes the proof.
\end{proof}

\begin{proof}[Proof of Proposition \ref{thm:analog_proposition1}]
Recall that
\begin{equation*}
\gwp(\lambda, \nu) = \inf_\pi \int c(\bx_1,\by_1,\bx_2,\by_2) ~ d \pi (\bx_1,\by_1) d \pi (\bx_2,\by_2).
\end{equation*}
We bound the above quantity by constructing an explicit transportation plan.  This is done as follows.  In the top level, fix a set $\cQ \in \mathfrak{Q}^1$, and consider the discrepancy $|\lambda(\cQ)-\mu(\cQ)|$.  This is the amount of mass that needs to be moved in and out of $Q$.  The cost in this move is bounded above by $h(1)$.  Altogether, this contributes a total of
$$
\begin{aligned}
& h(1) \times \Big( \sum_{\cQ \in \mathfrak{Q}^1} | \lambda(\cQ)-\mu(\cQ)| \Big) \times \Big( 2 \lambda (\cS) + \sum_{\cQ \in \mathfrak{Q}^1} | \lambda(\cQ)-\mu(\cQ)| \Big) \\
\leq \, & 3 h(1) \times \Big( \sum_{\cQ \in \mathfrak{Q}^1} | \lambda(\cQ)-\mu(\cQ)| \Big) \times \lambda (\cS).
\end{aligned}
$$

Now we consider moving mass within each of the $\cQ \in \mathfrak{Q}^1$.  We do so iteratively.  At each level of the refinement, the cost is bounded above by $c(\delta^{t-1})$.  By a similar set of arguments, the contribution is bounded above by
$$
3 h(\delta^{t-1}) \times \Big( \sum_{\cQ \in \mathfrak{Q}^t} | \lambda(\cQ)-\mu(\cQ)| \Big) \times \lambda (\cS).
$$

Finally, we consider the finest level of the discretization at $\mathfrak{Q}^{k}$.  We have to move all the mass within each $\cQ \in \mathfrak{Q}^k$.  The cost of moving is bounded above by $c(\delta^{k})$.  Then
\begin{equation*}
\begin{aligned}
& \sum_{\cQ_i \in \mathfrak{Q}^k, \cQ_j \in \mathfrak{Q}^k} \int_{\bx_1,\by_1 \in \cQ_i, \bx_2,\by_2 \in \cQ_j} c(\bx_1,\by_1,\bx_2,\by_2) ~ d \pi (\bx_1,\by_1) d \pi (\bx_2,\by_2) \\
\leq ~& h (\delta^k) \times \Big ( \sum_{\cQ_i \in \mathfrak{Q}^k} \int_{\bx_1,\by_1 \in \cQ_i} d \pi (\bx_1,\by_1) \Big )^2 \leq h (\delta^k) \times \lambda(\cS)^2.
\end{aligned}
\end{equation*}
This completes the proof.
\end{proof}

\begin{proposition}[\cite{WB:19}] \label{prop3-wb19}
Let $\cS$ be a Borel subset of $\cX$. Let $k^*$ be such that $\mathcal{N}_{3-k^*+1}(\cS) < +\infty$. Let $\{\cS_k\}_{k=1}^{k^*}$ be a sequence of Borel subsets of $\cS$. Then there exists a dyadic partition $\{\mathfrak{Q}_i^k\}_{k=1}^{k^*}$ of $\cS$ with parameter $\delta = 1/3$ such that for $1 \leq k \leq k^*$, the number of sets in $\mathfrak{Q}^k$ intersecting $\cS_k$ is at most $\mathcal{N}_{3-(k+1)}(\cS_k)$.
\end{proposition}

\begin{proposition}[\cite{WB:19}] \label{prop4-wb19}
Let $\cS$ be a Borel subset of $\cX$. Then
\begin{equation*}
\mathbb{E}\left[\sum_{\cQ_i^k \in \mathfrak{Q}_i^k} |\mu(\cQ_i^k) - \hat{\mu}_t(\cQ_i^k)|\right] \leq 2(1-\mu(\cS)) + \sqrt{|\{i: \cQ_i^k \cap \cS \neq \emptyset\}|/n}.
\end{equation*}
\end{proposition}

\begin{proof}[Proof of Theorem \ref{thm:convergencerate}]
First, we notice that for $n < \frac{9}{4} (\frac{3}{\varepsilon'} )^{s}$, the third term on the RHS will be greater than 1. The bound follows from $\gwp (\mu,\nu) \leq \text{diam}(\cX) \leq 1$ for measures on $\cX$. We thus focus on the case $n \geq \frac{9}{4} (\frac{3}{\varepsilon'} )^s$.

Define the following parameters $\alpha := \frac{sp}{s-2p}$, $\ell := \lceil\frac{-\log \varepsilon'}{\log 3}\rceil$, and $k^* := \lfloor\frac{\log n}{s\log 3}\rfloor - 2$.
Let $k'$ denote either the maximum integer in $[\ell, k^*]$ satisfying $k' \leq \frac{p}{\alpha} \cdot \frac{\log n}{s\log 3}$, or $\ell$ if no such integer exists. Our assumptions guarantee that for all $k \geq \ell$:
\begin{equation*}
\mathcal{N}_{3^{-k}}(\mu, 3^{-\alpha k}) \leq 3^{ks}.
\end{equation*}
Consequently, for all $k \geq k'$, we can find a subset $\cT_k \subset \cX$ with mass at least $1 - 3^{-\alpha k'}$ and at most 1 such that:
\begin{equation*}
\mathcal{N}_{3^{-k}}(\cT_k) \leq 3^{ks}.
\end{equation*}

By applying \Cref{thm:analog_proposition1} to $\cT_k$, taking $\delta = 3^{-1}$, and by using the fact that $c(\delta) \leq \delta^{pq}$ for the $L^{p,q}$ distortion distance, we have
\begin{equation*}
\mathbb{E}[ \gwp(\mu, \hat{\mu}_n)] ~\leq~ 3 \lambda(\cT_k) \times \sum_{k=0}^{k^*-1} 3^{-(k-1)pq} \Big (\sum_{Q^k_i \in \mathcal{Q}^{k}} | \lambda(\cQ^k_i) - \mu(\cQ^k_i)| \Big) + 3^{-kpq} \lambda(\cT_k)^2.
\end{equation*}

Then, by applying \Cref{prop4-wb19} to each of the term $\sum_{Q \in \mathcal{Q}^{k}} | \lambda(\cQ) - \mu(\cQ)|$ on the RHS, and by noting that $1 - \mu(\cT_k) \leq 3^{-\alpha k'}$ and $\mu(\cT_k) \leq 1$, we have 
\begin{equation}
\label{eq:intermediate-1}
\mathbb{E}[ \gwp(\mu, \hat{\mu}_n) ] ~\leq~ 3^{-k^*pq} + 3 \sum_{k=1}^{k^*-1} 3^{-(k-1)pq} \left(2\times 3^{-\alpha k} + \sqrt{|\{i: \cQ_i^k \cap \cT_k \neq \emptyset\}|/n} \right).
\end{equation}

Using \Cref{prop3-wb19} with $\cS_k = \cT_{k'}$ when $k < k'$ and $\cS_k = \cT_{k+1}$ when $k \geq k'$, we obtain a dyadic partition $\{\cQ^t\}_{1\leq t\leq k}$ where the number of $\cQ^k$ sets intersecting $\cT_k$ is bounded by $\mathcal{N}_{3^{-(k+1)}}(\cT_k)$. Replacing this bound to the RHS of \eqref{eq:intermediate-1}, we have
\begin{equation*}
\begin{aligned}
\mathbb{E}[ \gwp(\mu, \hat{\mu}_n) ] ~\leq~ & 3^{-k^*pq} + 3 \sum_{k=1}^{k'-1} 3^{-(k-1)pq}\sqrt{\frac{\mathcal{N}_{3^{-(k+1)}}(\cT_{k'})}{n}} \\
& \quad + 3 \sum_{k=k'}^{k^*} 3^{-(k-1)pq} \sqrt{\frac{\mathcal{N}_{3^{-(k+1)}}(\cT_{k+1})}{n}} 
+ 2 \cdot 3^{-\alpha k'}\sum_{k=1}^{k^*} 3^{-(k-1)pq}.
\end{aligned}
\end{equation*}

Since $\mathcal{N}_\varepsilon(\cT)$ increases monotonically as $\varepsilon$ decreases, we have the following for all $k \leq k' - 1$:
\begin{equation*}
\mathcal{N}_{3^{-(k+1)}}(\cT_{k'}) \leq \mathcal{N}_{3^{-k'}}(\cT_{k'}) \leq 3^{k's}.
\end{equation*}

By construction, for $k \geq k'$:
\begin{equation*}
\mathcal{N}_{3^{-(k+1)}}(\cT_{k+1}) \leq 3^{(k+1)s}.
\end{equation*}

Then, by combining these bounds together with the estimate $\sum_{k=1}^{k^*}3^{-(k-1)pq} \leq 3/2$ for $p,q \geq 1$ we have
\begin{equation*}
\begin{aligned}
    \mathbb{E}[ \gwp(\mu, \hat{\mu}_n) ] &\leq 3^{-k^*pq} + \frac{3}{2}\left(\dfrac{3 \times 3^{k's/2}}{\sqrt{n}} + 2 \times 3^{-\alpha k'} \right) + 3^{(1 + 2pq)} \sum_{k=k'}^{k^*}\dfrac{3^{(k+1)(s/2 - pq)}}{\sqrt{n}} \\
    &\leq 3^{-k^*pq} + \frac{3}{2}\left(\dfrac{3 \times 3^{k's/2}}{\sqrt{n}} + 2 \times 3^{-\alpha k'} \right) + \dfrac{3^{(1-k^*pq)}}{3^{s/2 - pq}-1}\sqrt{\dfrac{3^{(k^*+2)s}}{n}}.
\end{aligned}
\end{equation*}

From our choice of $k^*$, we have $3^{(k^*+2)s}\leq n$ and $3^{-k^*pq} \leq 3^{3pq}n^{-pq/s}$.  Moreover, the choice of $k'$ implies $\alpha k' > p\frac{\log n}{s\log 3} - 3\alpha$, which yields $3^{-\alpha k'} < 3^{3\alpha} n^{-p/s}$.  Therefore:

\begin{equation*}
    \mathbb{E}[ \gwp(\mu, \hat{\mu}_n) ] \leq \left( 3^{3pq} + \dfrac{3^{3pq}}{3^{s/2 - pq} - 1} + 3^{3\alpha + 1} \right) n^{-pq/s} + \frac{3 \times 3^{k's/2}}{2\sqrt{n}}.
\end{equation*}

Given that $k'$ satisfies $sk' \leq \max\{s\ell, (\frac{p}{\alpha}\frac{\log n}{\log 3} )\}$ we have
\begin{equation*}
    3^{k's/2} \leq 3^{\ell s/2} + n^{p/2\alpha} \leq 3^{\left(-\frac{\log \varepsilon'}{\log 3} + 1 \right) (s/2)} + n^{1/2}n^{-p/s} = \left(\frac{3}{\varepsilon'} \right)^{s/2} + n^{1/2}n^{-p/s}.
\end{equation*}

Finally, this yields our bound:
\begin{equation*}
    \mathbb{E}[ \gwp(\mu, \hat{\mu}_n) ] \leq \left( 3^{3pq} + \dfrac{3^{3pq}}{3^{s/2 - pq} - 1} + 3^{3\alpha + 1} \right) n^{-pq/s} + \frac{3n^{-p/s}}{2} + \frac{3}{2} \left(\frac{3}{\varepsilon'} \right)^{s/2}  n^{-1/2}.
\end{equation*}
\end{proof}

\section{Concluding Remarks}

In this paper, we introduce a sequence of optimization instances \eqref{hierarchy moment for cst} that relax the Gromov-Wasserstein problem \eqref{eq:gw}.  These optimization instances extend the moment-SOS hierarchy for the discretized GW problem in \eqref{eq:schmudgen_dgw} (and in \cite{tran:25}) in the sense that they are valid for general probability distributions in the source and target, and not only discrete ones.  The sequence of optimization instances \eqref{hierarchy moment for cst} satisfies a number of properties:  First, the discretized instances of these problems exactly recover the moment-SOS hierarchy \eqref{eq:schmudgen_dgw}.  Second, they provide a genuine hierarchy in the sense that the objective value increases and converges to the objective value of the Gromov-Wasserstein problem \eqref{eq:gw}.  Third, these optimization instances specify a pseudo-metric over the collection of metric measure spaces.  Separately from these properties, we also establish a statistical-consistency result that arises from sampling the source and target distributions uniformly at random.

\subsection{Future directions}

We pose some questions that stem from this work.

{\bf SOS for general optimization problems over distributions.}  The first question is whether there is a principled extension of the moment-SOS hierarchy to general optimization problems over probability distributions satisfying suitable conditions -- for instance, in settings where the objective and the constraints depend on the distributions (as optimization variables) in a polynomial way.  The GW problem \eqref{eq:gw} we study is one such example whereby the objective is quadratic, and the constraints are linear. 

{\bf PSD-ness for distributions.}  The second question is whether the continuous analog we identify should be considered the ``correct'' or the ``canonical'' extension.  One of our biggest conceptual difficulties in writing this paper is in formulating the constraint \eqref{cond: PSD for cts}.  Stated simply, what is the natural continuous analog of the PSD constraint that makes sense for distributions?  In a sense, our struggle stems from the following basic statement from linear algebra:  There are two equivalent ways to define a positive semidefinite matrix.  Let $X \in \mathbb{R}^{n \times n}$ be symmetric.  Then one has
\begin{equation*}
\bv^T X \bv \geq 0 \text{ for all } \bv \in \mathbb{R}^{n} \quad \Leftrightarrow \quad X = \textstyle \sum_{i} \lambda_i \bu_i \bu_i^T, \lambda_i \geq 0.
\end{equation*}
How do these definitions extend to linear operators?  The LHS characterization seems reasonably straightforward in which one replaces $\bv$ with functions, and vector-matrix sums with integrals -- this is precisely the approach we adopt in this paper.  The RHS characterization, unfortunately, does not lend itself to a natural extension.   To this end, we wish to point out ideas that appear to be superficially relevant:  namely, that compact, positive operators over Hilbert spaces admit a (unique) positive square-root \cite{Bol:99}.  However, we did not pursue this direction in this paper.

\bibliographystyle{alpha}
\bibliography{gwsoscts}

\end{document}